\documentclass[12pt,a4paper]{article}
\usepackage{}
\usepackage{mathrsfs}
\usepackage{indentfirst}
\usepackage{amsmath,amssymb,amsfonts,amsthm,graphics}
\usepackage{makeidx}
\usepackage{amsmath,amssymb,amsfonts,amsthm,graphics,graphicx}
\usepackage{makeidx}
\usepackage{color}
\usepackage{CJK}

\newtheorem{thm}{Theorem}[section]

\newtheorem{lem}{Lemma}[section]

\theoremstyle{definition}
\newtheorem{defn}{Definition}[section]

\def\-{\mbox{--}}

\newtheorem{claim}{Claim}[section]
\newtheorem{cor}{Corollary}[section]
\newtheorem{fact}{Fact}[section]

\newtheorem{remark}{Remark}[section]

\newtheorem{conj}{Conjecture}[section]

\def\qed{\hfill \nopagebreak\rule{5pt}{8pt}}
\def\pf{\noindent {\it Proof.} }

\title{Non-jumping Numbers for 5-Uniform Hypergraphs}

\author{Ran Gu, Xueliang Li, Zhongmei Qin, Yongtang Shi, Kang Yang \\
{\small  Center for Combinatorics and LPMC}\\
{\small Nankai University, Tianjin 300071, P.R. China}\\
{\small Email: guran323@163.com, lxl@nankai.edu.cn, qinzhongmei90@163.com,}\\
{\small shi@nankai.edu.cn, yangkang89@163.com} \\
}
\date{}
\begin{document}
\maketitle
\begin{abstract}
Let $\ell$ and $r$ be integers. A real number $\alpha \in [0,1)$ is a jump for $r$ if
for any $\varepsilon > 0$ and any integer $m,\ m \geq r$, any $r$-uniform graph with $n >
n_0(\varepsilon,m)$ vertices and at least $(\alpha+ \varepsilon)\binom{n}{r}$ edges
contains a subgraph with $m$ vertices and at least $(\alpha +c)\binom{m}{r}$ edges,
where $c=c(\alpha)$ is positive and does not depend on $\varepsilon$ and $m$. It follows from a theorem of Erd\H{o}s, Stone and Simonovits that every $\alpha \in [0,1)$ is a jump for $r=2$. Erd\H{o}s asked whether the same is true for $r \geq 3$. However, Frankl and R\"{o}dl gave a negative answer by showing that $1-\frac{1}{\ell^{r-1}}$ is not a jump for $r$ if $r \geq 3$ and $\ell >2r$. Peng gave more sequences of non-jumping numbers for $r=4$ and $r\geq 3$.  However, there are also a lot of unknowns on determining whether a number is a jump for $r \geq 3$. Following a similar approach as that of Frankl and R\"{o}dl, we give several sequences of non-jumping numbers for $r=5$, and extend one of the results to every $r \geq 5$, which generalize the above results.
\\[2mm]
\textbf{Keywords:} extremal problems in hypergraphs; Erd\H{o}s jumping constant conjecture;
Lagrangians of uniform graphs;
non-jumping numbers\\
\end{abstract}

\section{Introduction}
For a given finite set $V$ and a positive integer $r$, denote by $\binom{V}{r}$ the
family of all $r$-subsets of $V$. Let $G=(V(G),E(G))$ be a graph with {\it vertex set} $V(G)$ and {\it edge set} $E(G)$. We call $G$ an {\it $r$-uniform graph} if
$E(G)\subseteq \binom{V(G)}{r}$. An $r$-uniform graph $H$ is called a
{\it subgraph} of an $r$-uniform graph $G$ if $V(H)\subseteq V(G)$ and $E(H)\subseteq
E(G)$. Furthermore, $H$ is called an {\it induced subgraph} of $G$ if $E(H)=E(G)\cap
\binom{V(H)}{r}$.

Let $G$ be an $r$-uniform graph, we define the {\it density} of $G$ as
$\frac{|E(G)|}{\left|\binom{V(G)}{r}\right|}$, which is denoted by $d(G)$. Note that the density
of a complete $(\ell+1)$-partite graph with partition classes of size $m$ is greater
than $1-\frac{1}{\ell+1}$ (approaches $1-\frac{1}{\ell+1}$ when $m\rightarrow \infty$).
The density of a complete $r$-partite $r$-uniform graph with partition classes of size
$m$ is greater than $\frac{r!}{r^r}$ (approaches $\frac{r!}{r^r}$ when $m\rightarrow
\infty$).

In \cite{KNS}, Katona, Nemetz and Simonovits showed that, for any $r$-uniform graph
$G$, the average of densities of all induced subgraphs of $G$  with $m \geq r$
vertices is $d(G)$. From this result we know that there exists a subgraph of $G$
with $m$ vertices, whose density is at least $d(G)$. A natural question is: for a
constant $c>0$, whether there exists a subgraph of $G$ with $m$ vertices and density
at least $d(G)+c$ ? To be precise, the concept of ``jump" was introduced.

\begin{defn}
A real number $\alpha \in [0,1)$ is a jump for $r$ if there exists a constant $c>0$ such
that for any $\varepsilon>0$ and any integer $m, \ m \geq r$, any $r$-uniform graph with
$n> n_0(\varepsilon,m)$ vertices and density $\geq \alpha+\varepsilon$ contains a
subgraph with $m$ vertices and density $\geq \alpha +c$.
\end{defn}

Erd\H{o}s, Stone and Simonovits \cite{PM,PA} proved that every $\alpha \in [0,1)$ is a
jump for $r=2$. This result can be easily obtained from the following theorem.

\begin{thm}[\cite{PA}]
Suppose $\ell$ is a positive integer. For any $\varepsilon >0$ and any positive integer
$m$, there exists $n_0(m,\varepsilon)$ such that any graph $G$ on $n>n_0(m,\varepsilon)$
vertices with density $d(G) \geq 1-\frac{1}{\ell}+\varepsilon$ contains a copy of the
complete $(\ell +1)$-partite graph with partition classes of size $m$ (i.e., there
exists $\ell +1$ pairwise disjoint sets $V_1,\ldots,V_{\ell+1}$, each of them with size
$m$ such that $\{x,y\}$ is an edge whenever $x \in V_i$ and $y \in V_j$ for some $1 \leq
i <j \leq \ell+1$).
\end{thm}
\noindent Moreover, from the following theorem, Erd\H{o}s showed that for $r \geq
3$, every $\alpha \in [0,\frac{r!}{r^r})$ is a jump.

\begin{thm}[\cite{PE}]
For any $\varepsilon >0$ and any positive integer $m$, there exists $n_0(\varepsilon,m)$
such that any $r$-uniform graph $G$ on $n>n_0(\varepsilon,m)$ vertices with density
$d(G) \geq \varepsilon$ contains a copy of the complete $r$-partite $r$-uniform graph
with partition classes of size $m$ (i.e., there exist $r$ pairwise disjoint subsets
$V_1,\ldots, V_r,$ each of cardinality $m$ such that $\{x_1,x_2,\ldots,x_r\}$ is an edge
whenever $x_i \in V_i, 1\leq i \leq r$).
\end{thm}

\noindent Furthermore, Erd\H{o}s proposed the following conjecture on jumping
constant.
\begin{conj}
Every $\alpha \in [0,1)$ is a jump for every integer $r \geq 2$.
\end{conj}

Unfortunately,  Frankl and R\"{o}dl \cite{FR} disproved this conjecture by showing
the following result.

\begin{thm}[\cite{FR}]\label{thmFR}
Suppose $r \geq 3$ and $\ell > 2r$ are integers, then $1-\frac{1}{\ell^{r-1}}$ is not a jump for $r$.
\end{thm}

Using the approach developed by Frankl and R\"{o}dl in \cite{FR}, some other non-jump
numbers were given. However, for $r \geq 3$, there are still a lot of unknowns on determining whether a given
number is a jump. The breakthrough of the question is that Baber and Talbot \cite{BT} used Razborov's flag algebra method to show that jumps exist for $r=3$ in the interval $[2/9,1)$. These are the first examples of jumps for any $r \ge 3$ in the interval $[r!/r^r,1)$. They showed that for $r=3$ every $\alpha \in [0.2299, 0.2316)$ is a jump. Note that 0.2299 is very close to $\frac{2}{9}$, a well-known open question of Erd\H{o}s was raised.
\begin{center}{\it whether $\frac{r!}{r^r}$ \it is a jump for $r \geq 3$
and what is the smallest non-jump?}
\end{center}
In \cite{FPRT}, another question was raised:
\begin{center}{\it whether there is an interval of
non-jumps for some $r \geq 3$ ?}
\end{center}
Both questions seem to be very challenging. Regarding the first question, in
\cite{FPRT}, it was shown that $\frac{5r!}{2r^r}$ is a non-jump for $r \geq 3$  and
it is the smallest known non-jump until now. Some efforts were made in finding more
non-jumps for some $r \geq 3$. For $r = 3$, one more infinite sequence of non-jumps
(converging to 1)  was given in \cite{FPRT}. And for $r = 4$, several infinite
sequences of non-jumps (converging to 1) were found in \cite{P1,P2,P4,P5}. Every
non-jump in the above papers was extended to many sequences of non-jumps (still
converging to 1) in \cite{P3,P6,PZ}. Besides, in \cite{P7}, Peng found an infinite
sequence of non-jumps for $r = 3$ converging to $\frac{7}{12}$.

If a number $\alpha$ is a jump, then there exists a constant $c>0$ such that every
number in $[\alpha, \alpha + c)$ is a jump. As a direct result, we have that if
there is a set of non-jumping numbers whose limits form an interval (a number $a$ is
a limit of a set $A$ if there is a sequence $\{a_n\}_{n=1}^{\infty}, a_n \in A$ such
that $lim_{n\rightarrow\infty} a_n=a$), then every number in this interval is not a
jump. It is still an open problem whether such a ``dense enough" set of non-jumping
numbers exists or not.

In this paper, we intend to find more non-jumping numbers in addition to the known
non-jumping numbers given in \cite{FPRT,P1,P2,P3,P4,P5,P6,P7,PZ}.
Our approach is still based on the approach developed by Frankl and R\"{o}dl in
\cite{FR}. We first consider the case $r = 5$ and find a sequence of non-jumping
numbers. In Section $3$, we prove the following result.

\begin{thm}\label{thm1}
Let $\ell \geq 2$ be an integer. Then $1-\frac{5}{\ell^3}+\frac{4}{\ell^4}$ is not a jump for $r=5$.
\end{thm}

In \cite{P7}, Peng gave the following result: for positive integers $p \ge r \ge 3$, if $\alpha \cdot \frac{r!}{r^r}$ is a non-jump for $r$, then $\alpha \cdot \frac{p!}{p^p}$ is a non-jump for $p$. Combining with the Theorem \ref{thm1} for $\ell=5$ and $r=5$,
we have the following corollary directly.

\begin{cor}
Let $p \geq 5$, $\frac{151p!}{6p^p}$ is not a jump for $p$.
\end{cor}

We can use exactly the same reasoning to obtain more non-jumping numbers for $p>5$ by using Theorem \ref{thm1} together with Peng's result even with $\ell$ different from 5. We list it as the following corollary.

\begin{cor}
Let $p > 5$, $\frac{625}{24}(1-\frac{5}{\ell^3}+\frac{4}{\ell^4})\cdot \frac{p!}{p^p}$ is not a jump for $p$.
\end{cor}

Since in \cite{FPRT}, it was shown that $\frac{5r!}{2r^r}$ is a non-jumping number for
$r \geq 3$. In \cite{P3}, it was shown that for integers $r \geq 3$ and $p,\  3 \leq p
\leq r,\ (1-\frac{1}{p^{p-1}})\frac{p^p}{p!}\frac{r!}{r^r}$ is not a jump for $r$. In
particular, $\frac{12}{125}$ (take $r=5$ in $\frac{5r!}{2r^r}$), $\frac{96}{625}$ (take
$p=3$ and $r=5$ in $(1-\frac{1}{p^{p-1}})\frac{p^p}{p!}\frac{r!}{r^r}$) and
$\frac{252}{625}$ (take $p=4$ and $r=5$ in
$(1-\frac{1}{p^{p-1}})\frac{p^p}{p!}\frac{r!}{r^r}$) are non-jumping numbers for $r=5$.
In Section $4$, we will go back to the case of $r=5$ and prove the following result.
\begin{thm}\label{thm3}
Let $\ell \geq 2$ and $q \geq 1$ be integers. Then for $r=5$, we have

(a) If $q=1$ or $q \geq 2\ell^2+2\ell$, then $1-\frac{10}{\ell q}+\frac{35}{\ell^2
q^2}-\frac{50}{\ell^3 q^3}+\frac{4}{\ell^4 q^4}+\frac{10}{\ell q^4}-\frac{35}{\ell^2
q^4}+\frac{45}{\ell^3 q^4}$ is not a jump.

(b) If $q=1$ or $q \geq 10\ell^3$, then $1-\frac{10}{\ell q}+\frac{35}{\ell^2 q^2}
-\frac{50}{\ell^3 q^3}+\frac{10}{\ell q^4}-\frac{35}{\ell^2 q^4}+\frac{50}{\ell^3
q^4}-\frac{1}{\ell^4 q^4}$ is not a jump.

(c) $1-\frac{2}{q}+\frac{7}{5q^2}-\frac{2}{5q^3}+\frac{12}{125q^4}$ is not a jump.

(d) $1-\frac{2}{q}+\frac{7}{5q^2}-\frac{2}{5q^3}+\frac{96}{625q^4}$ is not a jump.

(e) If $q=1$ or $q \geq 3$, then
$1-\frac{2}{q}+\frac{7}{5q^2}-\frac{2}{5q^3}+\frac{252}{625q^4}$ is not a jump.
\end{thm}

When $q=1$, (a) reduces to Theorem \ref{thm1} for $r=5$, (b) reduces to Theorem
\ref{thmFR} for $r=5$, (c) shows that $\frac{12}{125}$ is not a jump for $r=5$, (d)
shows that $\frac{96}{625}$ is not a jump for $r=5$, and (e) shows that
$\frac{252}{625}$ is not a jump for $r=5$.

\section{Lagrangians and other tools}
In this section, we introduce the definition of Lagrangian of an $r$-uniform graph
and some other tools to be applied in the approach.

We first describe a definition of the Lagrangian of an $r$-uniform graph, which is a helpful
tool in the approach. More studies of Lagrangians were given in \cite{FF,FR,MS,T}.

\begin{defn}
For an $r$-uniform graph $G$ with vertex set $\{1,2,\ldots,m\}$, edge set $E(G)$ and a
vector $\vec{x}=\{x_1,\ldots,x_m\} \in R^m$, define
$$\lambda(G,\vec{x})=\sum_{\{i_1,\ldots,i_r\}\in E(G)}x_{i_1}x_{i_2}\cdots x_{i_r}.$$
$x_i$ is called the weight of vertex $i$.
\end{defn}

\begin{defn}
Let $S=\{\vec{x}=(x_1,x_2,\ldots,x_m):\sum_{i=1}^m x_i=1, \ x_i \geq 0$ for $
i=1,2,\ldots,m\}$. The Lagrangian of $G$, denoted by $\lambda(G)$, is defined as
$$\lambda(G)= max \{\lambda(G,\vec{x}):\vec{x} \in S\}.$$
\end{defn}

A vector $\vec{x}$ is called an optimal vector for $\lambda(G)$ if $\lambda(G,\vec{x})=\lambda(G)$.

We note that if $G$ is a subgraph of an $r$-uniform graph $H$, then for any vector
$\vec{x}$ in $S$, $\lambda(G,\vec{x}) \leq \lambda(H,\vec{x})$. The following fact is obtained directly.

\begin{fact}\label{fact1}
Let $G$ be a subgraph of an $r$-uniform graph $H$. Then
$$\lambda(G) \leq \lambda(H).$$
\end{fact}

For an $r$-uniform graph $G$ and $i \in V(G)$ we define $G_i$ to be the $(r-1)$-uniform
graph on $V-\{i\}$ with edge set $E(G_i)$ given by $e \in E(G_i)$ if and only if $e\cup
\{i\} \in E(G)$.

We call two vertices $i, j$ of an $r$-uniform graph $G$ {\it equivalent} if for all $f
\in \binom{V(G)-\{i,j\}}{r-1}$, $f \in E(G_i)$ if and only if $f \in E(G_j)$.

The following lemma given in \cite{P4} will be useful when calculating Lagrangians of some
certain hypergraphs.

\begin{lem}[\cite{P4}]\label{lem1}
Suppose $G$ is an $r$-uniform graph on vertices $\{1,2,\ldots,m\}$.
If vertices $i_1, i_2, \ldots, i_t$ are pairwise equivalent, then there exists an optimal
vector $\vec{y}=(y_1, y_2, \ldots, y_m)$ for $\lambda(G)$ such that $y_{i_1}=y_{i_2}=\cdots=y_{i_t}$.
\end{lem}

We also note that for an $r$-uniform graph $G$ with $m$ vertices, if we take
$\vec{x}=(x_1,x_2,\ldots,x_m)$, where each $x_i = \frac{1}{m}$, then for any fixed $\varepsilon >0$,  $$\lambda(G) \geq
\lambda(G,\vec{x})=\frac{|E(G)|}{m^r} \geq \frac{d(G)}{r!}-\varepsilon$$ for $m \geq
m'(\varepsilon)$.

On the other hand, we introduce a blow-up of an $r$-uniform graph $G$ which  allows
us to construct an $r$-uniform graph with a large number of vertices and density close
to $r!\lambda(G)$.

\begin{defn}
Let $G$ be an $r$-uniform graph with $V(G)=\{1,2,\ldots,m\}$ and
$\vec{\textbf{\emph{n}}}=(n_1,\ldots,n_m)$ be a positive integer vector. Define the
$\vec{\textbf{\emph{n}}}$ blow-up of $G$, $\vec{\textbf{\emph{n}}}\otimes G$ to be the
$m$-partite $r$-uniform graph with vertex set $V_1\cup \cdots \cup V_m, \ |V_i|=n_i, \ 1
\leq i \leq m$, and edge set $E(\vec{\textbf{\emph{n}}}\otimes
G)=\{\{v_{i_1},v_{i_2},\ldots,v_{i_r}\}:\ v_{i_k} \in V_{i_k} $ for $ 1 \leq k \leq r,\
\{i_1,i_2,\ldots,i_r\} \in E(G)\}$.
\end{defn}

In addition, we make the following easy remark given in \cite{P1}.

\begin{remark}[\cite{P1}]\label{rem1}
Let $G$ be an $r$-uniform graph with $m$ vertices and $\vec{y}=(y_1, y_2, \ldots, y_m)$
be an optimal vector for $\lambda(G)$. Then for any $\varepsilon>0$, there exists an
integer $n_1(\varepsilon)$, such that for any integer $n \geq n_1(\varepsilon)$,
\begin{equation}\label{equ1}
d((\lfloor ny_1\rfloor, \lfloor ny_2\rfloor,\ldots, \lfloor ny_m\rfloor)\otimes G)\geq
r! \lambda(G)-\varepsilon.
\end{equation}
\end{remark}

Let us also state a fact relating the Lagrangian of an $r$-uniform graph to the
Lagrangian of its blow-up used in \cite{FR} (\cite{FPRT,P1,P2,P4} as well).

\begin{fact}[\cite{FR}]\label{fact2}
If $n \geq 1$ and $\vec{\textbf{\emph{n}}}=(n,n,\ldots,n)$, then $\lambda
(\vec{\textbf{\emph{n}}}\otimes G)=\lambda(G)$ holds for every $r$-uniform graph $G$.
\end{fact}

We consider now the following definition.

\begin{defn}
For $\alpha \in [0,1)$ and a family $\mathcal {F}$ of $r$-uniform graphs, we say that
$\alpha$ is a threshold for $\mathcal {F}$ if for any $\varepsilon > 0$ there exists an
$n_0=n_0(\varepsilon)$ such that any $r$-uniform graph $G$ with $d(G) \geq
\alpha+\varepsilon$ and $|V(G)|>n_0$ contains some member of $\mathcal {F}$ as a
subgraph. We denote this fact by $\alpha \rightarrow \mathcal {F}$.
\end{defn}

The following lemma proved in \cite{FR} gives a necessary and sufficient condition for a
number $\alpha$ to be a jump.

\begin{lem}[\cite{FR}]\label{lem2}
The following two properties are equivalent.

1. $\alpha$ is a jump for $r$.

2. $\alpha \rightarrow \mathcal {F}$ for some finite family $\mathcal {F}$ of
$r$-uniform graphs satisfying $\lambda(F) > \frac{\alpha}{r!}$ for all $F \in \mathcal
{F}$.
\end{lem}

\begin{lem}[\cite{FR}]\label{lem3}
For any $\sigma \geq 0$ and any integer $k \geq r$, there exists $t_0(k,\sigma)$ such
that for every $t> t_0(k,\sigma)$, there exists an $r$-uniform graph $A$ satisfying:

1. $|V(A)|=t$.

2. $|E(A)| \geq \sigma t^{r-1}$.

3. For all $V_0 \subset V(A), r \leq |V_0| \leq k$ we have $|E(A) \cap \binom{V_0}{r}| \leq |V_0|-r+1$.
\end{lem}

We sketch the approach in proving Theorems \ref{thm1} and \ref{thm3}  as
follows (similar to the proof in \cite{P1,P2,P4}): Let $\alpha$ be the non-jumping
number described in those theorems. Assuming that $\alpha$ is a jump, we will derive a
contradiction by the following two steps.

\emph{Step 1}: Construct an $r$-uniform graph (in Theorems \ref{thm1}, \ref{thm3}, $r=5$)
with the Lagrangian close to but slightly smaller than $\frac{\alpha}{r!}$, then use
Lemma \ref{lem3} to add an $r$-uniform graph with a large enough number of edges but
spare enough (see properties 2 and 3 in Lemma \ref{lem3}) and obtain an $r$-uniform
graph with the Lagrangian  $\geq \frac{\alpha}{r!}+\varepsilon$ for some positive
$\varepsilon$. Then we ``blow up" this $r$-uniform graph to an new $r$-uniform graph, say
$H$, with a large enough number of vertices and density $>\alpha +
\frac{\varepsilon}{2}$ (see Remark \ref{rem1}). By Lemma
\ref{lem2}, if $\alpha$ is a jump then  $\alpha$ is a threshold for some finite family $\mathcal {F}$ of $r$-uniform
graphs with Lagrangian $> \frac{\alpha}{r!}$. So $H$ must contain some member of
$\mathcal {F}$ as a subgraph.

\emph{Step 2}: We show that any subgraph of $H$ with the number of vertices no more
than $\rm max \{|V(F)|,F \in \mathcal {F}\}$ has Lagrangian $ \leq \frac{\alpha}{r!}$
and derive a contradiction.

\section{Proof of Theorem
\ref{thm1}} In this section, we focus on $r=5$ and give a proof of Theorem
\ref{thm1}.

Let $\ell\geq 2$ and $\alpha=1-\frac{5}{\ell^3}+\frac{4}{\ell^4}$. Let $t$ be a
large enough integer given later. We first define a $5$-uniform hypergraph
$G(\ell,t)$ on $\ell$ pairwise disjoint sets $V_1, V_2,\ldots, V_\ell$, each of
cardinality $t$ whose density is close to $\alpha$ when $t$ is large enough. The
edge set of $G(\ell,t)$ consists of the following five kinds of hyperedges. When $\ell=2,3,4$, some of them may be vacant.

\begin{description}
\item[(1)] all $5$-subsets taking exactly one vertex from
each of $V_i$, $V_j$, $V_k$, $V_h$, $V_s$ ($1\leq i<j<k<h<s\leq \ell$);

\item[(2)] all $5$-subsets taking two vertices from $V_i$ and one vertex from each of $V_j$, $V_k$,
$V_h$ ($1\leq i\leq \ell$, $1\leq j<k<h\leq \ell$, $j,k,h\neq i$);

\item[(3)] all $5$-subsets taking two vertices from each of $V_i$, $V_j$ and one vertex from $V_k$ ($1\leq
i<j\leq \ell$, $1\leq k\leq \ell$, $k\neq i, j$);

\item[(4)] all $5$-subsets taking three
vertices from $V_i$, and one vertex from  each of $V_j$, $V_k$ ($1\leq i\leq \ell$,
$1\leq j<k\leq \ell$, $j,k\neq i$);

\item[(5)] all $5$-subsets taking three vertices from $V_i$
and two vertices from  $V_j$ ($1\leq i\leq \ell$, $1\leq j\leq \ell$, $j\neq i$).
\end{description}

Note that
\begin{align*}
|E(G(\ell,t ))|=&\binom{\ell}{5}t^5+\ell
\binom{\ell-1}{3}\binom{t}{2}t^3+\binom{\ell}{2}(\ell-2) \binom{t}{2}\binom{t}{2}t+
\ell\binom{\ell-1}{2}\binom{t}{3}t^2
\\[2mm]
&+\ell(\ell-1)
\binom{t}{3}\binom{t}{2}=\frac{\alpha}{120}\ell^5t^5
-c_0(\ell)t^4+o(t^4),
\end{align*}
where $c_0(\ell)$ is positive (we omit giving the precise calculation here). It is
easy to verify that the density of $G(\ell,t )$ is close to $\alpha$ if $t$ is large
enough. Corresponding to the $\ell t$ vertices of $G(\ell, t )$, we take the vector
$\vec x = (x_1, \ldots, x_{\ell t} )$, where $x_i =\frac{1}{\ell t}$ for each $i$,
$1\leq i \leq \ell t$, then
$$\lambda(G(\ell,t ))\geq \lambda(G(\ell,t ),\vec x)=\frac{|E(G(\ell,t ))|}{(\ell t)^5}
=\frac{\alpha}{120}-\frac{c_0(\ell)}{\ell^5t}+o(\frac{1}{t}),$$ which is close to
$\frac{\alpha}{120}$ when $t$ is large enough. We will use Lemma \ref{lem3} to add a
$5$-uniform graph to $G(\ell,t)$ so that the Lagrangian of the resulting graph is
$>\frac{\alpha}{120}+ \varepsilon(t)$ for some $\varepsilon(t)> 0$. Suppose that
$\alpha$  is a jump for $r=5$. According to Lemma \ref{lem2}, there exists a finite collection
$\mathcal{F}$ of $5$-uniform graphs satisfying:

 i) $\lambda(F ) >\frac{\alpha}{120}$
for all $F \in \mathcal{F}$, and

ii) $\alpha$ is a threshold for $\mathcal{F}$.

Set ${k_0} = ma{x_{_{_{F \in {\cal F}}}}}|V(F)|$ and $\sigma_0 = 2c_0(\ell)$. Let $r =
5$ and $t_0(k_0, \sigma_0)$ be given as in Lemma \ref{lem3}. Take an
integer $t >t_0$ and a $5$-uniform hypergraph $A(k_0,\sigma_0,t)$ satisfying the three
conditions in Lemma \ref{lem3}  with $V(A(k_0,\sigma_0,t))=V_1$. The $5$-uniform
hypergraph $H(\ell, t)$ is obtained by adding the hyperedges of $A(k_0,\sigma_0,t)$ to the $5$-uniform
hypegraph $G(\ell, t)$. For sufficiently large $t$, we have
$$\lambda(H(\ell, t))\geq \frac{|E(H(\ell, t))|}{(\ell t)^5} \geq \frac{|E(G(\ell, t))|
+\sigma_0t^4}{(\ell t)^5}\geq\frac{\alpha}{120}+\frac{c_0(\ell)}{2\ell^5t}.$$

Now suppose $\vec y = (y_1, y_2, \ldots, y_{\ell t} )$ is an optimal vector of $\lambda(H(\ell,
 t ))$. Let $ \varepsilon=\frac{30c_0(\ell )}{\ell^5t}$ and $n >n_1(\varepsilon)$ as in Remark
\ref{rem1}. Then the $5$-uniform graph $S_n = (\lfloor ny_1\rfloor, \ldots, \lfloor
ny_{\ell t}\rfloor)\otimes H(\ell,  t )$ has density not less than $\alpha+\varepsilon$.
Since $\alpha$ is a threshold for $\mathcal{F}$, some member $F$ of $\mathcal{F}$ is a
subgraph of $S_n$ for $n\geq max\{n_0(\varepsilon),n_1(\varepsilon)\}$. For such $F\in
\mathcal{F}$, there exists a subgraph $M$ of $H(\ell,  t )$ with $|V (M)|\leq |V (F
)|\leq k_0$, such that $F \subset \vec {\textbf{\emph{n}}} \otimes M$. By Facts \ref{fact1} and \ref{fact2}
we have
\begin{equation}\label{equ2}
\lambda(F )\leq \lambda( \vec {\textbf{\emph{n}}}\otimes M)= \lambda(M).
\end{equation}

\begin{lem}\label{lem4}
Let $M$ be any subgraph of $H(\ell ,t )$ with $|V (M)|\leq k_0$. Then
$$\lambda(M )\leq \frac{\alpha}{120}$$
holds.
\end{lem}
Applying Lemma \ref{lem4} to (\ref{equ2}), we have $\lambda(F)\leq \frac{\alpha}{120}$,
which contradicts our choice of $F$, i.e., contradicts the fact that $\lambda(F )
>\frac{\alpha}{120}$ for all $ F \in \mathcal{F}$.

{\em Proof of  Lemma \ref{lem4}. }By Fact \ref{fact1}, we may assume that $M$ is an
induced subgraph of $H(\ell ,t )$. Let $\widetilde{M}$ be an induced graph with maximum Lagrangian. Thus it is sufficient to show $\lambda(\widetilde{M})\leq \frac{\alpha}{120}$. Let $U_i =V (\widetilde{M})\cap V_i$. Define $M_1 =(U_1, E(\widetilde{M})\cap \binom{U_1}{5})$, i.e., the subgraph of $\widetilde{M}$ induced on $U_1$. In view of Fact \ref{fact1} and the assumption, it is enough to show Lemma \ref{lem4} for the case $E (M_1) \neq\emptyset$.
We assume $|V (M_1)| = 4 + d$ with $d$ a positive integer. By Lemma \ref{lem3},
$M_1$ has at most $d$ edges. Let $V (M_1) = \{v_1, v_2,\ldots, v_{4+d}\}$ and $\vec \xi
= (x_1,x_2,\ldots,x_{4+d})$ be an optimal vector for $\lambda(M_1)$ where $x_i$ is the
weight of vertex $v_i$. We may assume $x_1\geq x_2\geq \cdots\geq x_{4+d}$. The
following claim was proved in \cite{FR} ( Claim 4.4 there).
\begin{claim} \label{cl1}
$\sum\nolimits_{\left\{ {{v_i},{v_j},{v_k},{v_h},{v_s}} \right\} \in E\left( {{M_1}} \right)}
{{x_{{v_i}}}{x_{{v_j}}}{x_{{v_k}}}{x_{{v_h}}}{x_{{v_s}}}}  \leq \sum\limits_{5 \le i \le 4 + d}
{{x_1}{x_2}{x_3}{x_4}{x_i}}. $
\end{claim}

By Claim \ref{cl1} and the assumption, we may have $E(M_1) = \{\{v_1, v_2, v_3, v_4, v_i\}: \ 5\leq
i \leq 4+d\}$. Since $v_1, v_2, v_3, v_4$ are equivalent, in view of Lemma \ref{lem1},
we may assume that $x_1 =x_2 = x_3=x_4\buildrel {\rm def} \over = \rho$. For each $i$,
let $a_i$ be the sum of the weights of vertices of $U_i$. Notice that
\begin{equation*}
   \begin{cases}
   \sum\limits_{i = 1}^{\ell} {{a_i}}  = 1, \\
    a_i\geq 0,\ 1 \leq i \leq \ell \\
    0\leq\rho\leq\frac{a_1}{4}.
   \end{cases}
  \end{equation*}
Considering different types of edges in $M'$ and according to the definition of the
Lagrangian, we have
\begin{align*}
\lambda(\widetilde{M}) \leq& \sum\limits_{1 \le i < j < k < h < s \le \ell}
{{a_i}{a_j}{a_k}{a_h}{a_s}}  +\frac{1}{2}\sum\limits_{\scriptstyle2 \le i \le \ell;1 \le
j < k < h \le \ell;\hfill\atop
\scriptstyle j,k,h \ne i} {{a_i}^2{a_j}{a_k}{a_h}}   \\
\null  &\null+ \left(\sum\limits_{2 \le j < k < h \le \ell} {{a_j}{a_k}{a_h}}
\right)\left[\frac{1}{2}{({a_1} - 4\rho )^2} + 4\rho ({a_1} - 4\rho ) + 6{\rho
^2}\right]\\
\null  &\null  + \frac{1}{2}\left(\sum\limits_{\scriptstyle2 \le j \le \ell;2 \le k
\le \ell;\hfill\atop \scriptstyle k \ne j} {a_j^2{a_k}}\right)
\left[\frac{1}{2}{({a_1} - 4\rho )^2} + 4\rho ({a_1} - 4\rho ) + 6{\rho ^2}\right]
\\
\null  &\null  + \frac{1}{4}\sum\limits_{\scriptstyle2 \le i < j \le \ell;1 \le k
\le \ell;\hfill\atop \scriptstyle k \ne i,j}
{a_i^2a_j^2{a_k}}+\frac{1}{6}\sum\limits_{\scriptstyle2 \le i \le \ell;1 \le j < k
\le \ell;\hfill\atop
\scriptstyle j,k \ne i} {a_i^3{a_j}{a_k}}+ {\rho ^4}({a_1} - 4\rho )\\
\null  &\null  + \left(\sum\limits_{2 \le j < k \le \ell} {{a_j}{a_k}}\right)
\left[\frac{1}{6}{({a_1} - 4\rho )^3} + 2\rho
{({a_1} - 4\rho )^2} + 6{\rho ^2}({a_1} - 4\rho ) + 4{\rho ^3}\right]\\
\null  &\null+\frac{1}{{12}}\sum\limits_{\scriptstyle2 \le i \le \ell;2 \le j \le
\ell;\hfill\atop \scriptstyle j \ne i} {a_i^3a_j^2} +\frac{1}{6}\left(\sum\limits_{2
\le i \le \ell} {a_i^3}\right) \left[\frac{1}{2}{({a_1}
- 4\rho )^2} + 4\rho ({a_1} - 4\rho ) + 6{\rho ^2}\right] \\
\null  &\null+\frac{1}{2}\left(\sum\limits_{2 \le j \le \ell} {a_j^2}\right)
\left[\frac{1}{6}{({a_1} - 4\rho )^3} + 2\rho
{({a_1} - 4\rho )^2} + 6{\rho ^2}({a_1} - 4\rho ) + 4{\rho ^3}\right]\\
\null  =&\sum\limits_{1 \le i < j < k < h < s \le \ell} {{a_i}{a_j}{a_k}{a_h}{a_s}}
+ \frac{1}{2}\sum\limits_{\scriptstyle1 \le i \le \ell;1 \le j < k < h \le
\ell;\hfill\atop
\scriptstyle j,k,h \ne i} {{a_i}^2{a_j}{a_k}{a_h}} \\
\null  & + \frac{1}{4}\sum\limits_{\scriptstyle1 \le i < j \le \ell;1 \le k \le \ell;\hfill\atop
\scriptstyle k \ne i,j} {a_i^2a_j^2{a_k}}+\frac{1}{6}\sum\limits_{\scriptstyle1 \le i \le \ell;1 \le j < k
\le \ell;\hfill\atop
\scriptstyle j,k \ne i} {{a_i}^3{a_j}{a_k}}\null+\frac{1}{{12}}\sum\limits_{\scriptstyle1 \le i \le \ell;1
\le j \le \ell;\hfill\atop
\scriptstyle j \ne i} {{a_i}^3{a_j}^2}   \\
\null  & - 2{\rho ^2}\left(\sum\limits_{2 \le j < k < h \le \ell}
{{a_j}{a_k}{a_h}}\right)- {\rho ^2} \left(\sum\limits_{\scriptstyle2 \le j \le
\ell;2 \le k \le \ell;\hfill\atop \scriptstyle k \ne j} {a_j^2{a_k}}\right)  -
\frac{1}{3}{\rho ^2}\left(\sum\limits_{2 \le i \le \ell} {a_i^3} \right)\\
\null  &- {a_1}{\rho ^2}\left(\sum\limits_{2 \le j \le \ell} {a_j^2}\right) -
2{a_1}{\rho ^2}\left(\sum\limits_{2 \le j < k \le \ell} {{a_j}{a_k}}
\right)+\frac{4}{3}{\rho ^3} \left(\sum\limits_{2 \le j < k \le \ell} {{a_j}{a_k}}
\right)\\
\null  & + \frac{2}{3}{\rho ^3}\left(\sum\limits_{2 \le j \le \ell}
{a_j^2}\right) +
{\rho ^4}({a_1} - 4\rho )\\
\null  =&\sum\limits_{1 \le i < j < k < h < s \le \ell}
{{a_i}{a_j}{a_k}{a_h}{a_s}}  + \frac{1}{2}\sum\limits_{\scriptstyle1 \le i \le
\ell;1 \le j < k < h \le \ell;\hfill\atop
\scriptstyle j,k,h \ne i} {{a_i}^2{a_j}{a_k}{a_h}}
\end{align*}
\begin{align*}
\null  &+ \frac{1}{4}\sum\limits_{\scriptstyle1 \le i < j \le \ell;1 \le k \le \ell;\hfill\atop
\scriptstyle k \ne i,j} {a_i^2a_j^2{a_k}} +\frac{1}{6}\sum\limits_{\scriptstyle1 \le i \le \ell;1 \le j
 < k \le \ell;\hfill\atop
\scriptstyle j,k \ne i} {{a_i}^3{a_j}{a_k}}\null+\frac{1}{{12}}\sum\limits_{\scriptstyle1 \le i \le \ell;1
\le j \le \ell;\hfill\atop
\scriptstyle j \ne i} {{a_i}^3{a_j}^2}\\
\null  &- \frac{1}{3}{\rho ^2}\left(\sum\limits_{2 \le i \le \ell} {a_i} \right)^3-
{a_1}{\rho ^2}\left(\sum\limits_{2 \le i \le \ell} {a_i}\right)^2 +\frac{2}{3}{\rho
^3}\left(\sum\limits_{2 \le i \le \ell} {a_i}\right)^2  + {\rho ^4}\left({a_1} -
4\rho \right)\\
\null  =&\sum\limits_{1 \le i < j < k < h < s \le \ell} {{a_i}{a_j}{a_k}{a_h}{a_s}}
+ \frac{1}{2}\sum\limits_{\scriptstyle1 \le i\le \ell;1 \le j < k < h \le
\ell;\hfill\atop
\scriptstyle j,k,h \ne i} {a_i^2{a_j}{a_k}{a_h}}\\
\null  &+ \frac{1}{4}\sum\limits_{\scriptstyle1 \le i < j \le \ell;1 \le k \le \ell;\hfill\atop
\scriptstyle k \ne i,j} {a_i^2a_j^2{a_k}}+\frac{1}{6}\sum\limits_{\scriptstyle1 \le i \le \ell;1 \le j < k
\le \ell;\hfill\atop
\scriptstyle j,k \ne i} {{a_i}^3{a_j}{a_k}}\null+\frac{1}{{12}}\sum\limits_{\scriptstyle1 \le i \le \ell;1
\le j \le \ell;\hfill\atop
\scriptstyle j \ne i} {{a_i}^3{a_j}^2}\\
 \null  &+{\rho ^2}\left[{a_1}{\rho ^2} - 4{\rho ^3}+ \left(\frac{2}{3}\rho  - {a_1}\right)(1 - {a_1})^2 - \frac{1}{3}(1
 - {a_1})^3\right].
\end{align*}

Let
\begin{align}
&f(a_1 ,a_2,\ldots ,a_\ell,\rho)\nonumber\\
=&\sum\limits_{1 \le i < j < k < h < s \le \ell} {{a_i}{a_j}{a_k}{a_h}{a_s}}
+ \frac{1}{2}\sum\limits_{\scriptstyle1 \le i\le \ell;1 \le j < k < h \le
\ell;\hfill\atop
\scriptstyle j,k,h \ne i} {a_i^2{a_j}{a_k}{a_h}}\nonumber\\
\null  &+ \frac{1}{4}\sum\limits_{\scriptstyle1 \le i < j \le \ell;1 \le k \le \ell;\hfill\atop
\scriptstyle k \ne i,j} {a_i^2a_j^2{a_k}}+\frac{1}{6}\sum\limits_{\scriptstyle1 \le i \le \ell;1 \le j < k
\le \ell;\hfill\atop
\scriptstyle j,k \ne i} {{a_i}^3{a_j}{a_k}}\null+\frac{1}{{12}}\sum\limits_{\scriptstyle1 \le i \le \ell;1
\le j \le \ell;\hfill\atop
\scriptstyle j \ne i} {{a_i}^3{a_j}^2}\nonumber\\
 \null  &+{\rho ^2}\left[{a_1}{\rho ^2} - 4{\rho ^3}+ \left(\frac{2}{3}\rho  - {a_1}\right)(1 - {a_1})^2 - \frac{1}{3}(1
 - {a_1})^3\right].
\end{align}

Note that $$f(\frac{1}{\ell},\frac{1}{\ell},\ldots,\frac{1}{\ell},0)=\frac{\alpha}{120}.$$

Therefore, to show Lemma \ref{lem4}, we just need to show the following claim:

\begin{claim}\label{cla1}
$$f(a_1,a_2,\ldots,a_{\ell},\rho) \leq f(\frac{1}{\ell},\frac{1}{\ell},\ldots,\frac{1}{\ell},0)
=\frac{\alpha}{120}$$
holds under the constraints
\begin{equation*}
   \begin{cases}
   \sum\limits_{i = 1}^{\ell} {{a_i}}  = 1, \\
    a_i\geq 0, \ 1 \leq i \leq \ell\\
    0\leq\rho\leq\frac{a_1}{4}.
   \end{cases}
  \end{equation*}
\end{claim}

Define a fuction \begin{align*}
g(c_1 ,c_2,\ldots &,c_L) = \sum\limits_{1 \le i < j < k < h
< s \le L} {{c_i}{c_j}{c_k}{c_h}{c_s}}  + \frac{1}{2}\sum\limits_{\scriptstyle1 \le
i\le L;1 \le j < k < h \le L;\hfill\atop
\scriptstyle j,k,h \ne i} {c_i^2{c_j}{c_k}{c_h}}\\
+& \frac{1}{4}\sum\limits_{1 \le i < j \le L;1 \le k \le L;\hfill\atop \scriptstyle
k \ne i,j} {c_i^2c_j^2{c_k}}+\frac{1}{6}\sum\limits_{\scriptstyle1 \le i \le L;1 \le
j < k \le L;\hfill\atop \scriptstyle j,k \ne i}
{{c_i}^3{c_j}{c_k}}+\frac{1}{{12}}\sum\limits_{\scriptstyle1 \le i \le L;1 \le j \le
L;\hfill\atop \scriptstyle j \ne i} {{c_i}^3{c_j}^2},
\end{align*}
where $L \ge 2$ is an integer.

In order to prove Claim \ref{cla1}, we need to prove the following claim first:
\begin{claim} \label{cl2}
Let $c$ be a positive number. Suppose that
$\sum\limits_{i = 1}^L {{c_i}}
= c$ and each $c_i\geq 0$.
Then the function $g(c_1 ,c_2,\ldots,c_L)$
reaches the maximum $\frac{1}{120}(1-\frac{5}{L^3}+\frac{4}{L^4})c^5$ when $c_1 =
c_2 =\cdots= c_L = \frac{c}{L}$.
\end{claim}
\pf Since each term in function $g$ has degree $5$, we can assume that $c = 1$. Suppose
that $g$ reaches the maximum at $(c_1 , c_2 ,\ldots, c_L) $, we show that $c_1 = c_2
=\cdots= c_L = \frac{1}{L}$ must hold. If not, without loss of generality, assume that
$c_2 > c_1$, we will show that $g(c_1+\varepsilon, c_2-\varepsilon, c_3,\ldots,
c_L)-g(c_1, c_2,c_3,\ldots, c_L) > 0$ for small enough $\varepsilon > 0$ and derive a
contradiction. Notice that the summation of the terms in $g(c_1 , c_2 ,\ldots, c_L)$
containing $c_1 , c_2$ is
\begin{align*}
&(c_1 + c_2)\sum\limits_{3 \le i < j < k < h  \le L} {{c_i}{c_j}{c_k}{c_h}}  +c_1c_2 \sum\limits_{3 \le i < j < k \le L} {{c_i}{c_j}{c_k}}\\
&+\frac{1}{2}(c_1^2 + c_2^2)\sum\limits_{3 \le i < j < k \le L} {{c_i}{c_j}{c_k}}+\frac{1}{2}(c_1+ c_2)\sum\limits_{3 \le i \le L;3 \le j
<k\le L; j,k\neq i  } {{c_i^2}{c_j}{c_k}}\\
&+\frac{1}{2}(c_1^2c_2 + c_2^2c_1)\sum\limits_{3 \le i < j\le L} {{c_i}{c_j}}+\frac{1}{2}c_1 c_2\sum\limits_{3 \le i \le L;3 \le j \le L; j\neq i  }
 {{c_i^2}{c_j}}+\frac{1}{4}(c_1^2 c_2^2)\sum\limits_{3 \le i\le L} {{c_i}}\\
&+\frac{1}{4}(c_1^2c_2 + c_2^2c_1)\sum\limits_{3 \le i\le L} {{c_i^2}}+\frac{1}{4}(c_1^2+ c_2^2)\sum\limits_{3 \le i \le L;3 \le j \le L; j\neq i  }
 {{c_i^2}{c_j}}+\frac{1}{4}(c_1+ c_2)\sum\limits_{3 \le i < j\le L} {{c_i^2}{c_j^2}}\\
&+\frac{1}{6}(c_1^3+ c_2^3)\sum\limits_{3 \le i < j\le L} {{c_i}{c_j}}+\frac{1}{6}(c_1+ c_2)\sum\limits_{3 \le i \le L;3 \le j \le L; i\neq j }
{{c_i^3}{c_j}}+\frac{1}{6}(c_1^3c_2+ c_1c_2^3)\sum\limits_{3 \le i \le L} {{c_i}}\\
&+\frac{1}{6}c_1c_2\sum\limits_{3 \le i \le L} {{c_i^3}}+\frac{1}{12}(c_1^3+ c_2^3)\sum\limits_{3 \le i \le L} {{c_i^2}}+\frac{1}{12}(c_1^2+ c_2^2)
\sum\limits_{3 \le i \le L} {{c_i^3}}+\frac{1}{12}(c_1^3c_2^2+ c_1^2c_2^3)\\
=~&\frac{1}{24}(c_1+ c_2)[(\sum\limits_{3 \le i \le L} {{c_i}})^4-\sum\limits_{3 \le i \le L} {{c_i^4}}]+\frac{1}{12}(c_1+ c_2)^2
(\sum\limits_{3 \le i \le L} {{c_i}})^3\\
&+\frac{1}{12}(c_1+ c_2)^3(\sum\limits_{3 \le i \le L} {{c_i}})^2+\frac{1}{12}c_1c_2(2c_1^2+2c_2^2+3c_1c_2)\sum\limits_{3 \le i \le L} {{c_i}}+\frac{1}{12}
(c_1^3c_2^2+ c_1^2c_2^3)\\
=~&\frac{1}{24}(c_1+ c_2)(1-c_1-c_2)^4-\frac{1}{24}(c_1+ c_2)\sum\limits_{3 \le i \le L} {{c_i^4}}\\
&+\frac{1}{12}(c_1+ c_2)^2(1-c_1-c_2)^3+\frac{1}{12}(c_1+ c_2)^3(1-c_1-c_2)^2\\
&+\frac{1}{12}c_1c_2(2c_1^2+ 2c_2^2+3c_1c_2)(1-c_1-c_2)+\frac{1}{12}(c_1^3c_2^2+ c_1^2c_2^3).
\end{align*}
Therefore,
\begin{align*}
&\null g(c_1+\varepsilon, c_2-\varepsilon, c_3,\ldots, c_L)-g(c_1, c_2,c_3,\ldots, c_L)\\
=~&\frac{1}{12}(c_1+\varepsilon)(c_2-\varepsilon) [2(c_1+\varepsilon)^2+
2(c_2-\varepsilon)^2+3(c_1+\varepsilon)(c_2-\varepsilon)]
(1-c_1-c_2)\\
&+\frac{1}{12}(c_1+\varepsilon)^2
(c_2-\varepsilon)^2(c_1+ c_2)-\frac{1}{12}c_1c_2(2c_1^2+2c_2^2+3c_1c_2)(1-c_1-c_2)\\
&-\frac{1}{12}c_1^2c_2^2(c_1+c_2)\\
=~&\frac{1}{6}(c_2-c_1)(c_1^2+c_2^2+c_1c_2)(1-c_1-c_2)\varepsilon
+\frac{1}{6}c_1c_2(c_2-c_1)(c_1+ c_2)\varepsilon+o(\varepsilon)
>0.
\end{align*}
The last inequality is true since $c_1c_2$, $1-c_1-c_2$ cannot be equal to zero simultaneously (otherwise, $c_1=c_3=\cdots=c_L=0$ and $c_2=1$. Then $g(0,1,0,\ldots,0)=0<g(c/L, c/L, \ldots, c/L).$).
Therefore, $$g(c_1+\varepsilon, c_2-\varepsilon, c_3,\ldots, c_L)-g(c_1,
c_2,c_3,\ldots, c_L) > 0$$ for small enough $\varepsilon > 0$. This contradicts the
assumption that $g$ reaches the maximum at $(c_1 , c_2 ,\ldots, c_L) $.\qed\\

Since $0 \le \rho  \le \frac{{{a_1}}}{4},\ a_1 - 4\rho \geq 0,\ (1-a_1)^2 \geq 0$, then we have,
\begin{align*}
&{\rho ^2}\left[{a_1}{\rho ^2} - 4{\rho ^3} + \left(\frac{2}{3}\rho  - {a_1}\right){(1 - {a_1})^2} - \frac{1}{3}{(1 - {a_1})^3}\right]\\
\leq~&\rho ^2\left[\frac{a_1^3}{16 } - \frac{{{a_1^2}}}{4}{\rho } + \left(\frac{2}{3}\times\frac{{{a_1}}}{4}  - {a_1}\right){(1 - {a_1})^2}
- \frac{1}{3}{(1 - {a_1})^3}\right]\\
=~&\rho ^2\left[\frac{a_1^3}{16 } - \frac{{{a_1^2}}}{4}{\rho } - \left(\frac{{{a_1}}}{2}+\frac{1}{3}\right){(1 - {a_1})^2}\right]\\
=~&\rho ^2\left[\frac{1}{48}(-21{a_1^3}+32{a_1^2}+8a_1-16)-
\frac{1}{4}a_1^2\rho\right].
\end{align*}
Let $h(a_1)=-21{a_1^3}+32{a_1^2}+8a_1-16$, then, $h'(a_1)=-63{a_1^2}+64{a_1}+8$,
$h''(a_1)=-126{a_1}+64$. So $h'(a_1)$ increases when $0\leq a_1\leq\frac{32}{63}$,
$h'(a_1)$ decreases when $\frac{32}{63}\leq a_1\leq 1$. Hence, $h'(a_1)\geq
\min\{h'(0),h'(1)\}>0$, thus, $h(a_1)$ increases when $0\leq a_1\leq 1$. Note that
$h(0)<0$, $h(\frac{11}{15})<0$, $h(1)>0$, when $0\leq a_1\leq \frac{11}{15}$, we
have  ${\rho ^2}[{a_1}{\rho ^2} - 4{\rho ^3} + (\frac{2}{3}\rho  - {a_1}){(1 -
{a_1})^2} - \frac{1}{3}{(1 - {a_1})^3}]\leq 0$, by Claim \ref{cl2} and (3), we have
$f(a_1 ,a_2,\ldots ,a_\ell,\rho)\leq g(a_1 ,a_2,\ldots ,
a_\ell)\leq\frac{\alpha}{120}$. So Claim \ref{cla1} holds for $0 \leq a_1 \leq
\frac{11}{15}$. Therefore, we can assume that $\frac{11}{15}\leq a_1\leq 1$ and $\rho ^2\left[\frac{1}{48}(-21{a_1^3}+32{a_1^2}+8a_1-16)-
\frac{1}{4}a_1^2\rho\right]\ge 0$ (otherwise one could proceed as before for the case when $a_1\le 11/15$.). Since
the geometric mean is not greater than the arithmetic mean, we have,
\begin{align*}
&\rho ^2\left[\frac{a_1^3}{16 } - \frac{{{a_1^2}}}{4}{\rho } -
\left(\frac{{{a_1}}}{2}+\frac{1}{3}\right){(1 -
{a_1})^2}\right]\\
=~&\frac{64}{{{a_1^4}}}\left(\frac{{{a_1^2}\rho}}{8}\right)^2
\left[\frac{a_1^3}{16 } - \frac{{{a_1^2}}}{4}{\rho } - \left(\frac{{{a_1}}}{2}+\frac{1}{3}\right){(1 - {a_1})^2} \right]\\
\leq~&\frac{64}{{{a_1^4}}}\left[\frac{\frac{a_1^3}{16 } - \left(\frac{{{a_1}}}{2}+\frac{1}{3}\right){(1 - {a_1})^2}}{3}\right]^3\\
<~&\frac{64}{{{a_1^4}}}\left(\frac{a_1^3}{16\times3}\right)^3\leq \frac{1}{1728}.
\end{align*}

Let \begin{align*}
&k(a_1 ,a_2,\ldots ,a_\ell)\\
=&\sum\limits_{1 \le i < j < k < h < s \le \ell}
{{a_i}{a_j}{a_k}{a_h}{a_s}}+ \frac{1}{2}\sum\limits_{\scriptstyle1 \le i\le \ell;1
\le j < k < h \le
\ell;\hfill\atop \scriptstyle j,k,h \ne i} {a_i^2{a_j}{a_k}{a_h}}\\
&+ \frac{1}{4}\sum\limits_{\scriptstyle1 \le i < j \le \ell;1 \le k \le
\ell;\hfill\atop \scriptstyle k \ne i, j}
{a_i^2a_j^2{a_k}}+\frac{1}{6}\sum\limits_{\scriptstyle1 \le i \le \ell;1 \le j < k
\le \ell;\hfill\atop \scriptstyle j,k \ne i}
{{a_i^3}{a_j}{a_k}}+\frac{1}{{12}}\sum\limits_{\scriptstyle1 \le i \le \ell;1 \le j
\le \ell;\hfill\atop \scriptstyle j \ne i} {{a_i^3}{a_j^2}}+\frac{1}{1728}.
\end{align*}

Then combining with (3) we have
$$f(a_1 ,a_2,\ldots ,a_\ell,\rho)\leq k(a_1 ,a_2,\ldots ,a_\ell).$$
Therefore, to show Claim \ref{cla1}, it is sufficient to show the following claim:

\begin{claim} \label{cl3}
$$k(a_1 ,a_2,\ldots ,a_\ell)\leq\frac{\alpha}{120}$$ holds under the constraints $\sum\limits_{i = 1}^{\ell} {{a_i}}  = 1$,
$a_1\geq \frac{11}{15}$, and each $a_i\geq 0$.
\end{claim}

Define a function
\begin{align*}
r(a_2 ,a_3,\ldots ,a_\ell) =&\sum\limits_{2 \le  j < k < h < s \le \ell} {{a_j}{a_k}{a_h}{a_s}}+\frac{1}{2}\sum\limits_{\scriptstyle2 \le j\le \ell;2 \le  k < h \le \ell;\hfill\atop
\scriptstyle k,h \ne j} {a_j^2{a_k}{a_h}}+ \frac{1}{4}\sum\limits_{2 \le j < k \le \ell} {a_j^2a_k^2}\\
&+ \frac{1}{6}\sum\limits_{\scriptstyle2 \le j \le \ell;2 \le  k \le \ell;\hfill\atop
\scriptstyle k \ne j} {{a_j}^3{a_k}}.
\end{align*}
In order to prove Claim \ref{cl3}, we need to prove the following claim first:
\begin{claim} \label{cl4}
The function $r(a_2 ,a_3,\ldots ,a_\ell)$ reaches maximum $\frac{1}{24}(1-\frac{1}{(\ell-1)^3})c^4$ at
$a_2 = a_3 =\cdots= a_\ell = \frac{c}{\ell-1}$ under the constraints $\sum\limits_{i = 2}^{\ell} {{a_i}}  = c$, and each $a_i\geq 0$.
\end{claim}
{ \em Proof of  Claim \ref{cl4}. }Since $r(a_2 ,a_3,\ldots ,a_\ell)$ is a polynomial
with degree $4$ for each term, we just need to prove the claim for the case $c = 1$.
Suppose that $r$ reaches the maximum at $(c_2 , c_3 ,\ldots, c_\ell) $, we show that
$c_2 = c_3 =\cdots= c_\ell = \frac{1}{\ell-1}$. Otherwise, assume that $c_3 > c_2$, we
will show that $r(c_2+\varepsilon, c_3-\varepsilon, c_4, \ldots, c_\ell)-r(c_2 , c_3 ,\ldots,
c_\ell) > 0$ for small enough $\varepsilon > 0$ and derive a contradiction. Notice that
\begin{align*}
&r(c_2+\varepsilon, c_3-\varepsilon,c_4, \ldots, c_\ell)-r(c_2 , c_3 ,c_4,\ldots, c_\ell)\\
=~&[(c_2+\varepsilon)(c_3-\varepsilon)-c_2c_3]\sum\limits_{4 \le j < k \le \ell} {{c_j}{c_k}}\\
&+\frac{1}{2}[(c_2+\varepsilon)^2+ (c_3-\varepsilon)^2-c_2^2-c_3^2]\sum\limits_{4 \le j < k \le \ell} {{c_j}{c_k}}+\frac{1}{2}
[(c_2+\varepsilon) (c_3-\varepsilon)-c_2c_3]\sum\limits_{4 \le j \le \ell} {{c_j^2}}\\
&+\frac{1}{2}[(c_2+\varepsilon)^2(c_3-\varepsilon)+ (c_3-\varepsilon)^2(c_2+\varepsilon)-c_2^2c_3-c_2c_3^2]\sum\limits_{4 \le j  \le \ell} {{c_j}}\\
&+\frac{1}{4}[(c_2+\varepsilon)^2+ (c_3-\varepsilon)^2-c_2^2-c_3^2]\sum\limits_{4 \le j \le \ell} {{c_j^2}}+\frac{1}{4}[(c_2+\varepsilon)^2
(c_3-\varepsilon)^2-c_2^2c_3^2]\\
&+\frac{1}{6}[(c_2+\varepsilon)^3+ (c_3-\varepsilon)^3-c_2^3-c_3^3]\sum\limits_{4 \le j \le \ell} {{c_j}}\\
&+\frac{1}{6}[(c_2+\varepsilon)^3(c_3-\varepsilon)+
(c_3-\varepsilon)^3(c_2+\varepsilon)-c_2^3c_3-c_3^3c_2]\\
=~&\frac{1}{6}(c_3^3-c_2^3)\varepsilon+o(\varepsilon)>0,
\end{align*}
for small enough $\varepsilon > 0$ and get a contradiction.\qed\\

{ \em Proof of  Claim \ref{cl3}. }We will apply Claims \ref{cl2} and \ref{cl4}.
Separating the terms containing $a_1$ from the terms not containing $a_1$, we write
function $k(a_1 ,a_2,\ldots ,a_\ell)$ as follows:
\begin{align*}
&k(a_1 ,a_2,\ldots ,a_\ell) \\
=~&\sum\limits_{2 \le i < j < k < h < s \le \ell} {{a_i}{a_j}{a_k}{a_h}{a_s}}
+\frac{1}{2}\sum\limits_{\scriptstyle2 \le i\le \ell;2 \le j < k < h \le
\ell;\hfill\atop
\scriptstyle j,k,h \ne i} {a_i^2{a_j}{a_k}{a_h}}\\
+& \frac{1}{4}\sum\limits_{\scriptstyle2 \le i < j \le \ell;2 \le k \le \ell;\hfill\atop
\scriptstyle k \ne i, j} {a_i^2a_j^2{a_k}}+\frac{1}{6}\sum\limits_{\scriptstyle2 \le i \le \ell;2 \le j < k \le \ell;\hfill\atop
\scriptstyle j,k \ne i} {{a_i}^3{a_j}{a_k}}+\frac{1}{{12}}\sum\limits_{\scriptstyle2 \le i \le \ell;2 \le j \le \ell;\hfill\atop
\scriptstyle j \ne i} {{a_i}^3{a_j^2}}\\
+&a_1(\sum\limits_{2 \le  j < k < h < s \le \ell} {{a_j}{a_k}{a_h}{a_s}}+\frac{1}{2}\sum\limits_{\scriptstyle2 \le j\le \ell;2 \le  k < h \le \ell;\hfill\atop
\scriptstyle k,h \ne j} {a_j^2{a_k}{a_h}}+ \frac{1}{4}\sum\limits_{\scriptstyle2 \le j < k \le \ell} {a_j^2a_k^2}\\
+& \frac{1}{6}\sum\limits_{\scriptstyle2 \le j \le \ell;2 \le  k \le \ell;\hfill\atop
\scriptstyle k \ne j} {{a_j}^3{a_k}})+\frac{1}{{2}}a_1^2(\sum\limits_{2 \le  j < k < h \le \ell} {{a_j}{a_k}{a_h}})+\frac{1}{4}a_1^2(\sum\limits_{\scriptstyle2 \le j\le \ell;2 \le  k \le \ell;\hfill\atop
\scriptstyle k \ne j} {a_j^2{a_k}})\\
+&\frac{1}{6}a_1^3(\sum\limits_{\scriptstyle2 \le j <k \le \ell} {a_ja_k})+\frac{1}{12}a_1^3(\sum\limits_{2 \le j \le \ell} {a_j^2})+\frac{1}{12}a_1^2(\sum\limits_{2 \le j \le \ell} {a_j^3})+\frac{1}{1728}\\
=~&\sum\limits_{2 \le i < j < k < h < s \le \ell} {{a_i}{a_j}{a_k}{a_h}{a_s}}
+\frac{1}{2}\sum\limits_{\scriptstyle2 \le i\le \ell;2 \le j < k < h \le
\ell;\hfill\atop
\scriptstyle j,k,h \ne i} {a_i^2{a_j}{a_k}{a_h}}\\
+& \frac{1}{4}\sum\limits_{\scriptstyle2 \le i < j \le \ell;2 \le k \le \ell;\hfill\atop
\scriptstyle k \ne i, j} {a_i^2a_j^2{a_k}}+\frac{1}{6}\sum\limits_{\scriptstyle2 \le i \le \ell;2 \le j < k \le \ell;\hfill\atop
\scriptstyle j,k \ne i} {{a_i}^3{a_j}{a_k}}+\frac{1}{{12}}\sum\limits_{\scriptstyle2 \le i \le \ell;2 \le j \le \ell;\hfill\atop
\scriptstyle j \ne i} {{a_i}^3{a_j^2}}\\
+&a_1(\sum\limits_{2 \le  j < k < h < s \le \ell} {{a_j}{a_k}{a_h}{a_s}}+\frac{1}{2}\sum\limits_{\scriptstyle2 \le j\le \ell;2 \le  k < h \le \ell;\hfill\atop
\scriptstyle k,h \ne j} {a_j^2{a_k}{a_h}}+ \frac{1}{4}\sum\limits_{\scriptstyle2 \le j < k \le \ell} {a_j^2a_k^2}\\
+ &\frac{1}{6}\sum\limits_{\scriptstyle2 \le j \le \ell;2 \le  k \le \ell;\hfill\atop
\scriptstyle k \ne j} {{a_j}^3{a_k}})+\frac{1}{12}a_1^3(\sum\limits_{2 \le j \le \ell} {a_j})^2+\frac{1}{12}a_1^2(\sum\limits_{2 \le j \le \ell} {a_j})^3+\frac{1}{1728}.
\end{align*}
Applying Claim \ref{cl2} by taking $L = \ell-1 $ variables $a_2 , a_3 ,\ldots,a_\ell$
and $c = 1-a_1$, Claim \ref{cl4} and $\frac{1}{12}a_1^2(\sum\limits_{2 \le j \le \ell}
{a_j})^3+\frac{1}{12}a_1^3(\sum\limits_{2 \le j \le \ell}
{a_j})^2=\frac{1}{12}a_1^2(1-a_1)^3+\frac{1}{12}a_1^3(1-a_1)^2
=\frac{1}{12}a_1^2(1-a_1)^2$, we have
\begin{align*}
&k(a_1 ,a_2,\ldots ,a_\ell)\\
\le~&\frac{1}{120}\left[1-\frac{5}{(\ell-1)^3}
+\frac{4}{(\ell-1)^4}\right](1-a_1)^5\\
&+ \frac{1}{24}\left[1-\frac{1}{(\ell-1)^3}\right](1-a_1)^4a_1
+\frac{1}{12}a_1^2(1-a_1)^2+\frac{1}{1728}.
\end{align*}

Let \begin{align*}
w(a_1)=~&\frac{1}{120}\left[1-\frac{5}{(\ell-1)^3}
+\frac{4}{(\ell-1)^4}\right](1-a_1)^5\\
&+\frac{1}{24}\left[1-\frac{1}{(\ell-1)^3}\right](1-a_1)^4a_1
+\frac{1}{12}a_1^2(1-a_1)^2+\frac{1}{1728}.
\end{align*}
Therefore, to show Claim \ref{cl3}, we need to show the following claim:

\begin{claim} \label{cl5}
$$w(a_1)\leq\frac{\alpha}{120}$$ holds when $\frac{11}{15}\leq a_1\leq 1$.
\end{claim}
{ \em Proof.  } By a direct calculation,
$$w'(a_1)=\frac{1}{6}\left[\frac{1}{(\ell-1)^3}
-\frac{1}{(\ell-1)^4}\right](1-a_1)^4+\frac{1}{6(\ell-1)^3}
(1-a_1)^3a_1-\frac{1}{6}a_1^3(1-a_1),$$
$$w''(a_1)=\left[\frac{2}{3(\ell-1)^4}
-\frac{1}{2(\ell-1)^3}\right](1-a_1)^3
-\frac{1}{2(\ell-1)^3}(1-a_1)^2a_1+\frac{2}{3}a_1^3 -\frac{1}{2}a_1^2,$$
$$w^{(3)}(a_1)=\left[\frac{1}{(\ell-1)^3}
-\frac{2}{(\ell-1)^4}\right](1-a_1)^2 +\frac{1}{(\ell-1)^3}(1-a_1)a_1+2a_1^2 -a_1,$$
$$w^{(4)}(a_1)=\left[4
-\frac{4}{(\ell-1)^4}\right]a_1-1 -\frac{1}{(\ell-1)^3}+\frac{4}{(\ell-1)^4},$$ Note
that $w^{(4)}(a_1)>0$, when $\frac{11}{15}\leq a_1\leq 1$, so $w^{(3)}(a_1)$
increases when $\frac{11}{15}\leq a_1\leq 1$. By a direct calculation,
$w^{(3)}(\frac{11}{15})>0$,  so $w''(a_1)$ increases when $\frac{11}{15}\leq a_1\leq
1$. Since we have $w''(\frac{11}{15})<0$, $w''(1)>0$, thus, $w'(a_1)\leq
\max\{w'(\frac{11}{15}),w'(1)\}$.  By a direct calculation, $w'(\frac{11}{15})<0$,
$w'(1)=0$, so $w(a_1)$ is a decreasing function when $\frac{11}{15}\leq a_1\leq 1$.
When $\ell=2$, $w(\frac{11}{15})=\frac{1}{12}\times\frac{11^2\times4^2}{15^4}
+\frac{1}{1728}<\frac{1}{120}\times\frac{5}{8}=\frac{\alpha}{120}$. If $\ell\geq 3$,
since $1-\frac{5}{(\ell-1)^3} +\frac{4}{(\ell-1)^4}\geq1-\frac{5}{(2)^3}
+\frac{4}{(2)^4}$, then we have
$w(\frac{11}{15})=\frac{1}{120}[1-\frac{5}{(\ell-1)^3}
+\frac{4}{(\ell-1)^4}]-(1-\frac{4^5}{15^5})\times\frac{1}{120} (1-\frac{5}{(\ell)^3}
+\frac{4}{(\ell)^4})+\frac{1}{24}[1-\frac{1}{(\ell-1)^3}]\times\frac{11\times4^4}{15^5}+
\frac{1}{12}\times\frac{11^2\times4^2}{15^4}+\frac{1}{1728} \leq
\frac{1}{120}[1-\frac{5}{(\ell-1)^3}
+\frac{4}{(\ell-1)^4}]-(1-\frac{4^5}{15^5})\times\frac{1}{120} (1-\frac{5}{2^3}
+\frac{4}{2^4})+\frac{1}{24}\times\frac{11\times4^4}{15^5}+
\frac{1}{12}\times\frac{11^2\times4^2}{15^4}+\frac{1}{1728}
\leq\frac{1}{120}[1-\frac{5}{(\ell-1)^3} +\frac{4}{(\ell-1)^4}]$. So, $w(a_1)\leq
w(\frac{11}{15})\leq \frac{1}{120}[1-\frac{5}{(\ell-1)^3} +\frac{4}{(\ell-1)^4}]<
\frac{1}{120}[1-\frac{5}{\ell^3} +\frac{4}{\ell^4}]=\frac{\alpha}{120}$. This
completes the proof of Claim \ref{cl5}.\qed\\

Applying Claim \ref{cla1} to (3), we have
$$\lambda(\widetilde{M})\leq\frac{\alpha}{120}.$$ This completes the proof of Lemma \ref{lem4}. \qed

\section{Proof of Theorem \ref{thm3}}
In this section, we focus on $r=5$ and prove the following Theorem, which implies
Theorem \ref{thm3}.

\begin{thm}\label{thm4}
Let $\ell \geq 2,\ q \geq 1$ be integers. Let $N(\ell)$ be any of the five numbers given below.

\begin{equation}\label{equa1}
N(\ell)=\left\{ \begin{array}{ll}
\alpha(1-\frac{5}{\ell^3}+\frac{4}{\ell^4}), \ or\\
1-\frac{1}{\ell^4}, \ or\\
\frac{12}{125} \ (in \ this \ case, \ view \ \ell=5), \ or \\
\frac{96}{625} \ (in \ this \ case, \ view  \ \ell=5), \ or\\
\frac{252}{625} \ (in \ this \ case, \ view \ \ell=5).
\end{array} \right.
\end{equation}
Then
\begin{equation}\label{equa2}
N(\ell,q)=1-\frac{10}{\ell q}+\frac{35}{\ell^2 q^2}-\frac{50}{\ell^3 q^3}+\frac{10}{\ell q^4}
-\frac{35}{\ell^2 q^4}+\frac{50}{\ell^3 q^4}-\frac{1}{q^4}+\frac{N(\ell)}{q^4}
\end{equation}
is not a jump for 5 provided
\begin{equation}\label{equa3}
q=1 \ or \ \ell^3(1-N(\ell))(q^3+q^2+q+1)-10\ell^2(q^2+q+1)+35\ell(q+1)-50 \geq 0
\end{equation}
holds.
\end{thm}

Now let us explain why Theorem \ref{thm4} implies Theorem \ref{thm3}.

If $N(\ell)=\alpha$, then
\begin{align*}
&\ell^3(1-N(\ell))(q^3+q^2+q+1)-10\ell^2(q^2+q+1)+35\ell(q+1)-50\\
=~&\ell^3(\frac{5}{\ell^3}-\frac{4}{\ell^4})(q^3+q^2+q+1)-10\ell^2(q^2+q+1)+35\ell(q+1)-50\\
=~&\frac{1}{\ell}[(5\ell-4)q^3+(5\ell-10\ell^3-4)q^2+(5\ell-10\ell^3+35\ell^2-4)q\\
 &+(-45\ell-10\ell^3+35\ell^2-4)].
\end{align*}
Let
\begin{align*}
f_1(q)=~&\frac{1}{\ell}[(5\ell-4)q^3+(5\ell-10\ell^3-4)q^2+(5\ell-10\ell^3+35\ell^2-4)q\\
 &+(-45\ell-10\ell^3+35\ell^2-4)],
\end{align*}
then $f_1(q)$ is an increasing function of $q$ when $q \geq 2\ell^2+2\ell$ and $f_1(2\ell^2+2\ell)>0$.
Therefore, when $q \geq 2\ell^2+2\ell$, (\ref{equa3}) is satisfied. Applying Theorem
\ref{thm4}, we get Part (a) of Theorem \ref{thm3}.

If $N(\ell)=1-\frac{1}{\ell^4}$, then
\begin{align*}
&\ell^3(1-N(\ell))(q^3+q^2+q+1)-10\ell^2(q^2+q+1)+35\ell(q+1)-50\\
=~&\ell^3(\frac{1}{\ell^4})(q^3+q^2+q+1)-10\ell^2(q^2+q+1)+35\ell(q+1)-50\\
=~&\frac{1}{\ell}[q^3-(10\ell^3-1)q^2-(10\ell^3-35\ell^2-1)q+(1-10\ell^3+35\ell^2-50\ell)].
\end{align*}
Let
$$f_2(q)=\frac{1}{\ell}[q^3-(10\ell^3-1)q^2-(10\ell^3-35\ell^2-1)q+(1-10\ell^3+35\ell^2-50\ell)],$$
then $f_2(q)$ is an increasing function of $q$ when $q \geq 7\ell^3$ and $f_2(10\ell^3)>0$.
Therefore, when $q \geq 10\ell^3$, (\ref{equa3}) is satisfied. Applying Theorem
\ref{thm4}, we get Part (b) of Theorem \ref{thm3}.

If $\ell=5$ and $N(\ell)=\frac{12}{125}$, then
\begin{align*}
&\ell^3(1-N(\ell))(q^3+q^2+q+1)-10\ell^2(q^2+q+1)+35\ell(q+1)-50\\
=~&113q^3-137q^2+38q-12.
\end{align*}
Let
$$f_3(q)=113q^3-137q^2+38q-12,$$
then $f_3(q)$ is an increasing function of $q$ when $q \geq 1$ and $f_3(1)>0$. Therefore, (\ref{equa3}) is satisfied.
Applying Theorem \ref{thm4}, we get Part (c) of Theorem \ref{thm3}.

If $\ell=5$ and $N(\ell)=\frac{96}{625}$, then
\begin{align*}
&\ell^3(1-N(\ell))(q^3+q^2+q+1)-10\ell^2(q^2+q+1)+35\ell(q+1)-50\\
=~&\frac{1}{5}(529q^3-721q^2+154q-96).
\end{align*}
Let
$$f_4(q)=\frac{1}{5}(529q^3-721q^2+154q-96),$$
then $f_4(q)$ is an increasing function of $q$ when $q \geq 1$ and $f_4(2)>0$. Therefore, (\ref{equa3}) is satisfied.
Applying Theorem \ref{thm4}, we get Part (d) of Theorem \ref{thm3}.

If $\ell=5$ and $N(\ell)=\frac{252}{625}$, then
\begin{align*}
&\ell^3(1-N(\ell))(q^3+q^2+q+1)-10\ell^2(q^2+q+1)+35\ell(q+1)-50\\
=~&\frac{1}{5}(373q^3-877q^2-2q-252).
\end{align*}
Let
$$f_5(q)=\frac{1}{5}(373q^3-877q^2-2q-252),$$
then $f_5(q)$ is an increasing function of $q$ when $q \geq 2$ and $f_5(3)>0$. Therefore, when $q \geq 3$, (\ref{equa3}) is satisfied.
Applying Theorem \ref{thm4}, we get Part (e) of Theorem \ref{thm3}.

Now we give the proof of Theorem \ref{thm4}.

\emph{Proof of Theorem \ref{thm4}.} We will show that $N(\ell,q)$ is not
a jump for $5$. Let $t$ be a fixed large enough integer determined later. We first
define a $5$-uniform hypergraph $G(\ell,t)$ on $\ell$ pairwise disjoint sets
$V_1,\ldots, V_{\ell}$, each of them with size $t$ and the density of $G(\ell,t)$ is
close to $N(\ell)$ when $t$ is large enough. Each of the five choices of $N(\ell)$
corresponds to a construction.

1. If $N(\ell)=\alpha$, then $G(\ell,t)$ is defined in section $3$. Notice that
$$d(G(\ell,t))=\frac{\binom{\ell}{5}t^5+\ell
\binom{\ell-1}{3}\binom{t}{2}t^3+\binom{\ell}{2}(\ell-2) \binom{t}{2}\binom{t}{2}t+
\ell\binom{\ell-1}{2}\binom{t}{3}t^2
+\ell(\ell-1)
\binom{t}{3}\binom{t}{2}}{\binom{\ell t}{5}}$$
which is close to $\alpha$ if $t$ is large enough.

2. If $N(\ell)=1-\frac{1}{\ell^4}$, then $G(\ell,t)$ is defined on $\ell$ pairwise
disjoint sets $V_1,V_2,\ldots, V_{\ell}$, where $|V_i|=t$, and the edge set of
$G(\ell,t)$ is $\binom{\cup_{i=1}^{\ell}V_i}{5}-\cup_{i=1}^{\ell}\binom{V_i}{5}$. Notice
that
$$d(G(\ell,t))=\frac{\binom{\ell t}{5}-\ell \binom{t}{5}}{\binom{\ell t}{5}}$$
which is close to $1-\frac{1}{\ell^4}$ if $t$ is large enough.

3. If $N(5)=\frac{12}{125}$ (in this case, view $\ell=5$), then $G(5,t)$ is defined on
$5$ pairwise disjoint sets $V_1,V_2,V_3,V_4,V_5$, where $|V_i|=t$, and the edge set of
$G(5,t)$ consists of all 5-sets in the form of $\{\{a,b,c,v_4,v_5\}, $ where $ a \in
V_1, \ b \in V_2,\ c \in V_3 $ and $ v_4 \in V_4, \ v_5 \in V_5 \}$, or
$\{\{a,b,c,v_4,v_5\}, $ where $ \{a,b\} \in \binom{V_1}{2},\ c \in V_2 $ and $ v_4 \in
V_4, \ v_5 \in V_5 \}$, or $\{\{a,b,c,v_4,v_5\}, $ where $ \{a,b\} \in \binom{V_2}{2}, \
c \in V_3 $ and $ v_4 \in V_4, \ v_5 \in V_5 \}$, or $\{\{a,b,c,v_4,v_5\}, $ where $
\{a,b\} \in \binom{V_3}{2},\ c \in V_1 $ and $ v_4 \in V_4,\ v_5 \in V_5 \}$. Notice
that
$$d(G(5,t))=\frac{t^5+3\binom{t}{2}t^3}{\binom{5t}{5}}$$
which is close to $\frac{12}{125}$ if $t$ is large enough.

4. If $N(5)=\frac{96}{625}$ (in this case, view $\ell=5$), then $G(5,t)$ is defined on
5 pairwise disjoint sets $V_1,V_2,V_3,V_4,V_5$, where $|V_i|=t$, and the edge set of
$G(5,t)$ consists of all 5-sets in the form of $\{\{v_1,v_2,v_3,v_4,v_5\},$ where $
\{v_1,v_2,v_3\} \in \binom{\cup_{i=1}^{3}V_i}{3}-\cup_{i=1}^{3}\binom{V_i}{3}, $ and $
v_4 \in V_4, \ v_5 \in V_5\}$. Notice that
$$d(G(5,t))=\frac{(\binom{3t}{3}-3\binom{t}{3})t^2}{\binom{5t}{5}}$$
which is close to $\frac{96}{625}$ if $t$ is large enough.

5. If $N(5)=\frac{252}{625}$ (in this case, view $\ell=5$), then $G(5,t)$ is defined on
5 pairwise disjoint sets $V_1,V_2,V_3,V_4,V_5$, where $|V_i|=t$, and the edge set of
$G(5,t)$ consists of all 5-sets in the form of $\{\{v_1,v_2,v_3,v_4,v_5\},$ where $
\{v_1,v_2,v_3,v_4\} \in \binom{\cup_{i=1}^{4}V_i}{4}-\cup_{i=1}^{4}\binom{V_i}{4}, $ and
$ v_5 \in V_5\}$. Notice that
$$d(G(5,t))=\frac{(\binom{4t}{4}-4\binom{t}{4})t}{\binom{5t}{5}}$$
which is close to $\frac{252}{625}$ if $t$ is large enough.

We also note that
\begin{equation}\label{equa4}
\frac{|E(G(\ell,t))|+\frac{1}{12}\ell^4 t^4}{(\ell t)^5}\geq \frac{1}{120}(N(\ell)+\frac{1}{\ell^5 t})
\end{equation}
holds for $t \geq t_1$.

The 5-uniform graph $G(\ell,q,t)$ on $\ell q$ pairwise disjoint sets $V_i, \ 1 \leq i
\leq \ell q$, each of them with size $t$ is obtained as follows: for each $p,\ 0 \leq p
\leq q-1$, take a copy of $G(\ell,t)$ on the vertex set $\cup_{p \ell +1 \leq j \leq
(p+1) \ell}V_j$, then add all other edges in
the form of $\{\{v_{j_1},v_{j_2},v_{j_3},v_{j_4},v_{j_5}\}, $ where $ 1 \leq j_1 <j_2
<j_3 <j_4 <j_5 \leq \ell q $ and $ v_{j_k} \in V_{j_k} $ for $ 1 \leq k \leq 5$, and $\cup_{1\le k \le 5}v_{j_k} \nsubseteq \cup_{p \ell +1 \leq j \leq
(p+1) \ell}V_j\}$. We
will use Lemma \ref{lem3} to add a 5-uniform graph to $G(\ell,q,t)$ so that the
Lagrangian of the resulting graph is $> \frac{N(\ell,q)}{120}+\varepsilon(t)$ for some
$\varepsilon(t)>0$. The precise argument is given below.

Suppose that $N(\ell,q)$ is a jump for $r=5$. By Lemma \ref{lem2}, there exists a finite
collection $\mathcal {F}$ of 5-uniform graphs satisfying the following:

i) $\lambda(F)>\frac{N(\ell,q)}{120}$ for all $F \in \mathcal {F}$, and

ii) $N(\ell,q)$ is a threshold for $\mathcal {F}$.

Assume that $r=5$ and set $k_1= \max _{F \in \mathcal {F}}|V(F)|$ and
$\sigma_1=\frac{1}{12}\ell^4 q$. Let $t_0(k_1,\sigma_1)$ be given as in Lemma
\ref{lem3}. Fix an integer $t> \rm max(t_0,t_1)$, where $t_1$ is the number from
(\ref{equa4}).

Take a 5-uniform graph $A_{k_1,\sigma_1}(t)$ satisfying the conditions in Lemma
\ref{lem3} with $V(A_{k_1,\sigma_1}(t))=V_1$. The 5-uniform hypergraph $H(\ell,q,t)$ is
obtained by adding $A_{k_1,\sigma_1}(t)$ to the 5-uniform hypergraph $G(\ell,q,t)$. Now
we give a lower bound of $\lambda(H(\ell,q,t))$. Notice that,
$$\lambda(H(\ell,q,t))\geq \frac{|E(H(\ell,q,t))|}{(\ell qt)^5}.$$

In view of the construction of $H(\ell,q,t)$, we have
\begin{align*}
&\frac{|E(H(\ell,q,t))|}{(\ell qt)^5}\ge\frac{|E(G(\ell,q,t))|+\sigma_1 t^4}{(\ell qt)^5}\\
=~&\frac{q|E(G(\ell,t))|+\frac{1}{12}\ell^4 q t^4 + (\binom{\ell q}{5}-q \binom{\ell}{5})t^5}{(\ell qt)^5}\\
=~&\frac{q|E(G(\ell,t))|+\frac{1}{12}\ell^4 q t^4}{(\ell qt)^5} + \frac{1}{120}(1-\frac{10}{\ell q}+\frac{35}{\ell^2 q^2}
-\frac{50}{\ell^3 q^3}-\frac{1}{q^4}+\frac{10}{\ell q^4}-\frac{35}{\ell^2 q^4}+\frac{50}{\ell^3 q^4})\\
\mathrel{\mathop \geq^{(\ref{equa4})}}~&\frac{1}{120}(\frac{N(\ell)}{q^4}+\frac{1}{(\ell q)^5 t})+\frac{1}{120}(1-\frac{10}
{\ell q}+\frac{35}{\ell^2 q^2}-\frac{50}{\ell^3 q^3}-\frac{1}{q^4}+\frac{10}{\ell q^4}-\frac{35}{\ell^2 q^4}+\frac{50}{\ell^3 q^4})\\
\mathrel{\mathop =^{(\ref{equa2})}}~&\frac{1}{120}(N(\ell,q)+\frac{1}{(\ell q)^5
t}).
\end{align*}

Hence, we have
$$\lambda(H(\ell,q,t))\geq \frac{1}{120}(N(\ell,q)+\frac{1}{(\ell q)^5 t}). $$

Now suppose $\vec{y}=\{y_1,y_2,\ldots,y_{\ell qt}\}$ is an optimal vector of
$\lambda(H(\ell,q,t))$. Let $\varepsilon =\frac{1}{2(\ell q)^5 t}$ and
$n>n_1(\varepsilon)$ as in Remark \ref{rem1}. Then the 5-uniform graph $S_n=(\lfloor
ny_1\rfloor,\ldots,\lfloor ny_{\ell qt}\rfloor ) \otimes H(\ell,q,t)$ has density
larger than $N(\ell,q)+\varepsilon$. Since $N(\ell,q)$ is a threshold for $\mathcal
{F}$, some member $F$ of $\mathcal {F}$ is a subgraph of $S_n$ for $n\geq \rm
max\{n_0(\varepsilon),n_1(\varepsilon)\}$. For such $F \in \mathcal {F}$, there exists
a subgraph $M'$ of $H(\ell,q,t)$ with $|V(M')| \leq k_1$ so that $F \subset
\vec{\textbf{\emph{n}}}\otimes M' \subset \vec{\emph{\textbf{n}}}\otimes H(\ell,q,t)$.

Theorem \ref{thm4} will follow from the following lemma.

\begin{lem}\label{lem6}
Let $M'$ be any graph of $H(\ell,q,t)$ with $|V(M')| \leq k_1$. Then
\begin{equation}\label{equa5}
\lambda(M') \leq \frac{1}{120}N(\ell,q)
\end{equation}
holds.
\end{lem}

The proof of Lemma \ref{lem6} will be given as follows. We continue the proof of
Theorem \ref{thm4} by applying this Lemma. By Fact \ref{fact2} we have
$$\lambda(F) \leq \lambda(\vec{\textbf{\emph{n}}}\otimes M') = \lambda(M') \leq
\frac{1}{120}N(\ell,q)$$ which contradicts our choice of $F$, i.e., contradicts the
fact that $\lambda(F)>\frac{1}{120}N(\ell,q)$ for all $F \in \mathcal {F}$. This
completes the proof of Theorem \ref{thm4}.\qed\\

\emph{Proof of Lemma \ref{lem6}.} Let $M'$ be any subgraph of $H(\ell,q,t)$ with
$|V(M')|\leq k_1$ and $\vec{\xi}$ be an optimal vector for $\lambda(M')$. Define
$U_i=V(M')\cap V_i$ for $1 \leq i \leq \ell q$. Let $a_i$ be the sum of the weights in
$U_i, \ 1 \leq i \leq \ell q$, respectively. Note that $\sum_{i=1}^{\ell q}a_i=1$ and
$a_i \geq 0$ for each $i,\ 1 \leq i \leq \ell q$.

The proof  of Lemma \ref{lem6} is based on Lemma \ref{lem4}, Claims \ref{cla1},
\ref{cl2} and an estimation given in \cite{FPRT} and \cite{P3} on the summation of
the terms in $\lambda(M')$ corresponding to edges in $E(M')\cap
\binom{\cup_{i=1}^{\ell}V_i}{5}$, denoted by $\lambda(M'\cap\cup _{i=1}^{\ell}V_i)$.
For our purpose, we formulate Claim \ref{cla1} in Section 3, Lemma 4.2 in
\cite{FPRT} and Lemma 3.2 in \cite{P3} as follows.

\begin{lem}\label{lem7}
There exists a function $f$ such that
\begin{equation}\label{equa6}
\lambda(M'\cap\cup _{i=1}^{\ell}V_i) \leq f(a_1,a_2,\ldots,a_{\ell},\rho),
\end{equation}
where the function $f$ satisfies the following property:
\begin{equation}\label{equa7}
f(a_1,a_2,\ldots,a_{\ell},\rho)\leq
f(\frac{c}{\ell},\frac{c}{\ell},\ldots,\frac{c}{\ell},0)=\frac{1}{120}N(\ell)c^5
\end{equation}
under the constraints $\sum_{j=1}^{\ell}a_j=c$ and each $a_j \geq 0,\ 1 \leq j
\leq \ell$ for any positive constant $c$ , and for positive constant $\rho$, $0 \leq \rho \leq \frac{a_1}{4}$.
\end{lem}

In view of the construction of $H(\ell,q,t)$, for each $p,\  1\leq p \leq q-1$, the
structure of $M'$ restricted on the vertex set $\cup_{i=p \ell +1}^{(p+1)\ell}V_i$ is
similar to the structure of $M'$ restricted on the vertex set $\cup_{i=1}^{\ell}V_i$, but
there might be some other extra edges in $\binom{V_1}{5}$ for $M'$ restricted on the
vertex set $\cup_{i=1}^{\ell}V_i$. Therefore, for each $p, \ 1\leq p \leq q-1$, let the
summation of the terms in $\lambda(M')$ corresponding to edges in $E(M')\cap
\binom{\cup_{i=p \ell +1}^{(p+1)\ell}V_i}{4}$ denoted by $\lambda(M'\cap\cup _{i=p \ell +1}^{(p+1)\ell}V_i)$. For
our purpose, we formulate Claim \ref{cl2} in section 3, Lemma 4.2 in \cite{FPRT} and
Lemma 3.2 in \cite{P3} as follows.

\begin{lem}\label{lem8}
For every integer $p$, there exists a function $g$ such that
\begin{equation}\label{equa8}
\lambda(M'\cap\cup _{i=p \ell+1}^{(p+1)\ell}V_i) \leq g(a_{p \ell+1},a_{p \ell +2},\ldots, a_{(p+1)\ell}),
\end{equation}
where the function $g$ satisfies the following property:
\begin{equation}\label{equa9}
g(d_{p \ell+1},d_{p \ell +2},\ldots, d_{(p+1)\ell})\leq g(\frac{c}{\ell},\frac{c}{\ell},\ldots,\frac{c}{\ell})=\frac{1}{120}N(\ell)c^5
\end{equation}
under the constraints $\sum_{j=p \ell+1}^{(p+1)\ell}d_j=c$ and each $d_j \geq 0, \ p \ell+1 \leq j \leq (p+1)\ell$ for any positive constant $c$.
\end{lem}
Consequently,
\begin{align*}
\lambda(M')& \leq ~f(a_1,a_2,\ldots,a_{\ell},\rho)+\sum_{p=1}^{q-1}g(a_{p \ell+1},a_{p \ell +2},\ldots, a_{(p+1)\ell})\\
+&(\sum_{1 \leq i_1<i_2<i_3<i_4<i_5 \leq \ell q}a_{i_1}a_{i_2}a_{i_3}a_{i_4}a_{i_5}-\sum_{p=0}^{q-1}\sum_{p \ell+1 \leq i_1<i_2<i_3<i_4<i_5 \leq (p+1)\ell}a_{i_1}a_{i_2}a_{i_3}a_{i_4}a_{i_5}).
\end{align*}

Let \begin{align*}
&F(a_1,a_2,\ldots,a_{\ell q},\rho) \\
=~&f(a_1,a_2,\ldots,a_{\ell},\rho)+\sum_{p=1}^{q-1}g(a_{p \ell+1},a_{p \ell +2},\ldots, a_{(p+1)\ell})\\
&+(\sum_{1 \leq i_1<i_2<i_3<i_4<i_5 \leq \ell q}a_{i_1}a_{i_2}a_{i_3}a_{i_4}a_{i_5}-\sum_{p=0}^{q-1}\sum_{p \ell+1 \leq i_1<i_2<i_3<i_4<i_5 \leq (p+1)\ell}a_{i_1}a_{i_2}a_{i_3}a_{i_4}a_{i_5}).
\end{align*}

Note that
\begin{equation}\label{equa10}
F(\frac{1}{\ell q},\frac{1}{\ell q},\ldots,\frac{1}{\ell q},0)=\frac{N(\ell)}{120q^4}+\frac{\binom{\ell q}{5}-q\binom{\ell}{5}}{(\ell q)^5}=\frac{N(\ell,q)}{120}.
\end{equation}
Therefore, to show Lemma \ref{lem6}, we only need to show the following claim:

\begin{claim}\label{cl6}
\begin{equation}\label{equa11}
F(a_1,a_2,\ldots,a_{\ell q},\rho)\leq F(\frac{1}{\ell q},\frac{1}{\ell q},\ldots,\frac{1}{\ell q},0)
\end{equation}
holds under the constraints $\sum_{j=1}^{\ell q}a_j=1$ and each $a_j \geq 0, \ 1 \leq j \leq \ell q$ and $0 \leq \rho \leq \frac{a_1}{4}$.
\end{claim}

\pf Suppose the function $F$ reaches the maximum at $(a_1,a_2,\ldots,a_{\ell},\rho)$.
By applying Lemma \ref{lem7}, we claim that we can assume that $a_1=a_2=\cdots=a_{\ell}$
and $\rho=0$. Otherwise, let
$c_1=c_2=\cdots=c_{\ell}=\frac{\sum_{j=1}^{\ell}a_j}{\ell}$. Then
\begin{align*}
&F(c_1,c_2,\ldots,c_{\ell},a_{\ell+1},\ldots,a_{\ell q},0)-F(a_1,a_2,\ldots,a_{\ell},a_{\ell+1},\ldots,a_{\ell q},\rho)\\
=~&f(c_1,c_2,\ldots,c_{\ell},0)-f(a_1,a_2,\ldots,a_{\ell},\rho)\\
&+(\sum_{1\leq i<j<k<h\leq \ell}c_i c_j c_k c_h-\sum_{1\leq i<j<k<h\leq \ell}a_i a_j a_k a_h)(\sum_{s=\ell +1}^{\ell q}a_s)\\
&+(\sum_{1\leq i<j<k\leq \ell}c_i c_j c_k -\sum_{1\leq i<j<k\leq \ell}a_i a_j a_k )(\sum_{\ell +1 \leq h < s \leq \ell q}a_h a_s)\\
&+(\sum_{1\leq i<j\leq \ell}c_i c_j -\sum_{1\leq i<j\leq \ell}a_i a_j)(\sum_{\ell +1 \leq k < h < s \leq \ell q}a_k a_h a_s)
\geq 0
\end{align*}
holds by combining (\ref{equa7}), $\sum_{1\leq i<j<k<h\leq \ell}c_i c_j c_k c_h-\sum_{1\leq i<j<k<h\leq \ell}a_i a_j a_k a_h \geq 0$ ,\\$\sum_{1\leq i<j<k\leq \ell}c_i c_j c_k -\sum_{1\leq i<j<k\leq \ell}a_i a_j a_k \geq 0$ and $\sum_{1\leq i<j\leq \ell}c_i c_j -\sum_{1\leq i<j\leq \ell}a_i a_j \geq 0$. This implies that $a_1=a_2=\cdots=a_{\ell}$ and $\rho=0$ can be assumed. Similarly, by applying Lemma \ref{lem8}, for each $p, \ 1 \leq p \leq q-1$, we can assume that $a_{p \ell +1}=a_{p\ell +2}= \cdots =a_{(p+1)\ell}$. Set $b_{p+1}=a_{p \ell +1}=a_{p\ell +2}= \cdots =a_{(p+1)\ell}$ for each $0 \leq p \leq q-1$. Let
\begin{align*}
&R(b_1,b_2,\ldots,b_q)\\
=~&\frac{N(\ell)}{120}\sum_{p=1}^{q}\ell^5 b_p^5+\sum_{p=1}^{q}\binom{\ell}{4}b_p^4(1-\ell b_p)+\sum_{1\leq p_1 \leq q;1\leq p_2 \leq q;p_2\neq p_1}\binom{\ell}{3}\binom{\ell}{2}b_{p_1}^3 b_{p_2}^2\\
&+\sum_{1 \leq p_1 \leq q; 1 \leq p_2 < p_3 \leq q; p_2,p_3 \neq p_1}\binom{\ell}{3}\ell^2 b_{p_1}^3 b_{p_2}b_{p_3}+\sum_{1\leq p_1<p_2 \leq q; 1 \leq p_3 \leq q; p_3 \neq p_1,p_2}\binom{\ell}{2}^2 \ell b_{p_1}^2 b_{p_2}^2 b_{p_3}\\
&+\sum_{1 \leq p_1 \leq q; 1 \leq p_2<p_3<p_4 \leq q; p_2,p_3,p_4 \neq p_1}\binom{\ell}{2} \ell^3 b_{p_1}^2 b_{p_2}b_{p_3}b_{p_4}+\sum_{1 \leq p_1<p_2<p_3<p_4<p_5 \leq q}\ell^5 b_{p_1} b_{p_2} b_{p_3} b_{p_4} b_{p_5}.
\end{align*}

In view of Lemma \ref{lem7} and Lemma \ref{lem8}, we have
$$F(a_1,a_2,\ldots,a_{\ell q},\rho)\leq R(b_1,b_2,\ldots,b_q).$$

Note that
\begin{equation}\label{equa12}
R(\frac{1}{\ell q},\frac{1}{\ell q},\ldots,\frac{1}{\ell q})=F(\frac{1}{\ell q},\frac{1}{\ell q},\ldots,\frac{1}{\ell q},0)\mathrel{\mathop =^{(\ref{equa10})}}\frac{N(\ell,q)}{120}.
\end{equation}
Therefore, to show Claim \ref{cl6}, it is sufficient to show the following claim£»

\begin{claim}\label{cl7}
$$R(b_1,b_2,\ldots,b_q)\leq R(\frac{1}{\ell q},\frac{1}{\ell q},\ldots,\frac{1}{\ell q})$$
holds under the constraints
\begin{equation}\label{equa13}
\left\{ \begin{array}{ll}
\sum_{i=1}^q b_i=\frac{1}{\ell},\\
b_i \geq 0, \quad 1 \leq i \leq q.
\end{array} \right.
\end{equation}
\end{claim}

Suppose that function $R$ reaches the maximum at $(b_1,b_2,\ldots,b_q)$. We will apply
Claims \ref{cl8} and \ref{cl9} stated below.

\begin{claim}\label{cl8}
Let $i,\ j,\ 1 \leq i <j \leq q$ be a pair of integers and $\varepsilon$ be a real
number. Let $c_i=b_i + \varepsilon,\ c_j=b_j - \varepsilon$, and $c_k = b_k$ for $k \neq
i, j$. Let $(b_j-b_i)A(b_1, b_2,\ldots,b_q)$ and $B(b_1,b_2,\ldots,b_q)$ be the
coefficients of $\varepsilon$ and $\varepsilon^2$ in
$R(c_1,c_2,\ldots,c_q)-R(b_1,b_2,\ldots,b_q)$, respectively, i.e.,
\begin{align*}
&R(c_1,c_2,\ldots,c_q)-R(b_1,b_2,\ldots,b_q)\\
=&(b_j-b_i)A(b_1,
b_2,\ldots,b_q)\varepsilon+B(b_1,b_2,\ldots,b_q)\varepsilon^2+o(\varepsilon^2).
\end{align*}
If $b_i \neq b_j$, then
$$A(b_1, b_2,\ldots,b_q)+B(b_1,b_2,\ldots,b_q) \geq 0.$$
\end{claim}

\pf Without loss of generality, we take $i=1$ and $j=2$. By the definition of
the function $R(b_1,b_2,\ldots,b_q)$, we have
\begin{align*}
&R(b_1+\varepsilon,b_2-\varepsilon,\ldots,b_q)-R(b_1,b_2,\ldots,b_q)\\
=~&\frac{N(\ell)}{120}\ell^5[(b_1+\varepsilon)^5+(b_2-\varepsilon)^5-b_1^5-b_2^5]\\
&+\binom{\ell}{4}[(b_1+\varepsilon)^4(1-\ell b_1-\ell \varepsilon)+(b_2-\varepsilon)^4(1-\ell b_2+\ell \varepsilon)-b_1^4(1-\ell b_1)-b_2^4(1-\ell b_2)]\\
&+\binom{\ell}{3}\binom{\ell}{2}[(b_1+\varepsilon)^3+(b_2-\varepsilon)^3-b_1^3-b_2^3](\sum_{3 \leq p_1\leq q}b_{p_1}^2)\\
&+\binom{\ell}{3}\binom{\ell}{2}[(b_1+\varepsilon)^2+(b_2-\varepsilon)^2-b_1^2-b_2^2](\sum_{3 \leq p_1 \leq q}b_{p_1}^3)\\
&+\binom{\ell}{3}\binom{\ell}{2}[(b_1+\varepsilon)^3(b_2-\varepsilon)^2
+(b_2-\varepsilon)^3(b_1+\varepsilon)^2-b_1^3b_2^2-b_2^3b_1^2]\\
&+\binom{\ell}{3}\ell^2[(b_1+\varepsilon)^3+(b_2-\varepsilon)^3-b_1^3-b_2^3](\sum_{3 \leq p_1 < p_2\leq q}b_{p_1}b_{p_2})\\
&+\binom{\ell}{3}\ell^2[(b_1+\varepsilon)^3(b_2-\varepsilon)+(b_2-\varepsilon)^3(b_1+\varepsilon)-b_1^3b_2-b_2^3b_1]
(\sum_{3 \leq p_1\leq q}b_{p_1})\\
&+\binom{\ell}{3}\ell^2[(b_1+\varepsilon)(b_2-\varepsilon)-b_1b_2](\sum_{3\leq p_1 \leq q}b_{p_1}^3)\\
&+\binom{\ell}{2}^2\ell[(b_1+\varepsilon)^2+(b_2-\varepsilon)^2-b_1^2-b_2^2](\sum_{3\leq p_1 \leq q; 3 \leq p_2 \leq q; p_2 \neq p_1}b_{p_1}^2b_{p_2})\\
&+\binom{\ell}{2}^2\ell[(b_1+\varepsilon)^2(b_2-\varepsilon)^2-b_1^2b_2^2](\sum_{3 \leq p_1\leq q}b_{p_1})\\
&+\binom{\ell}{2}^2\ell[(b_1+\varepsilon)^2(b_2-\varepsilon)+(b_2-\varepsilon)^2(b_1+\varepsilon)-b_1^2b_2-b_2^2b_1]
(\sum_{3 \leq p_1\leq q}b_{p_1}^2)\\
&+\binom{\ell}{2}\ell^3[(b_1+\varepsilon)^2+(b_2-\varepsilon)^2-b_1^2-b_2^2](\sum_{3 \leq p_1 < p_2<p_3\leq q}b_{p_1}b_{p_2}b_{p_3})\\
&+\binom{\ell}{2}\ell^3[(b_1+\varepsilon)^2(b_2-\varepsilon)+(b_2-\varepsilon)^2(b_1+\varepsilon)-b_1^2b_2-b_2^2b_1](\sum_{3 \leq p_1 < p_2\leq q}b_{p_1}b_{p_2})\\
&+\binom{\ell}{2}\ell^3[(b_1+\varepsilon)(b_2-\varepsilon)-b_1b_2](\sum_{3\leq p_1 \leq q; 3 \leq p_2 \leq q; p_2 \neq p_1}b_{p_1}^2b_{p_2})\\
&+\ell^5[(b_1+\varepsilon)(b_2-\varepsilon)-b_1b_2](\sum_{3 \leq p_1 < p_2<p_3\leq q}b_{p_1}b_{p_2}b_{p_3}).
\end{align*}

By a direct calculation, we obtain that
\begin{align*}
&A(b_1, b_2,\ldots,b_q)+B(b_1,b_2,\ldots,b_q)\\
=~&-\frac{N(\ell)}{24}\ell^5(b_1+b_2)(b_1^2+b_2^2)+5\ell \binom{\ell}{4}(b_1+b_2)(b_1^2+b_2^2)-4\binom{\ell}{4}(b_1^2+b_2^2+b_1b_2)\\
&+2\binom{\ell}{3}\binom{\ell}{2}b_1b_2(b_1+b_2)+\binom{\ell}{3}\ell^2(b_1-b_2)^2(\sum_{3\leq p_1 \leq q}b_{p_1})+2\binom{\ell}{2}^2\ell b_1b_2(\sum_{3\leq p_1 \leq q}b_{p_1})\\
&+\frac{N(\ell)}{12}\ell^5(b_1^3+b_2^3)+\binom{\ell}{4}(6b_1^2+6b_2^2-10 \ell b_1^3 -10\ell b_2^3)+\binom{\ell}{3}\binom{\ell}{2}(b_1^3+b_2^3-3b_1b_2^2-3b_1^2b_2)\\
&-3\binom{\ell}{3}\ell^2(b_1-b_2)^2(\sum_{3\leq p_1 \leq q}b_{p_1})+\binom{\ell}{2}^2\ell(b_1^2+b_2^2-4b_1b_2)(\sum_{3\leq p_1 \leq q}b_{p_1})\\
=~&[2\ell \binom{\ell}{4}-2\binom{\ell}{3}\ell^2+\binom{\ell}{2}^2\ell]\frac{1}{\ell}(b_1-b_2)\\
&+[\frac{N(\ell)}{24}\ell^5-5\ell\binom{\ell}{4}+\binom{\ell}{3}\binom{\ell}{2}+2\binom{\ell}{3}\ell^2-\binom{\ell}{2}^2\ell]
(b_1+b_2)(b_1-b_2)^2\\
&\geq [2\ell \binom{\ell}{4}-2\binom{\ell}{3}\ell^2+\binom{\ell}{2}^2\ell](b_1+b_2)(b_1-b_2)^2\\
&+[\frac{N(\ell)}{24}\ell^5-5\ell\binom{\ell}{4}+\binom{\ell}{3}\binom{\ell}{2}+2\binom{\ell}{3}\ell^2-\binom{\ell}{2}^2\ell]
(b_1+b_2)(b_1-b_2)^2\\
=~&[\frac{N(\ell)}{24}\ell^5-3\ell\binom{\ell}{4}+\binom{\ell}{3}\binom{\ell}{2}](b_1+b_2)(b_1-b_2)^2\\
=~&\left\{ \begin{array}{ll}
(\frac{5}{12}\ell^4-\frac{23}{24}\ell^3+\frac{3}{8}\ell^2+\frac{1}{6}\ell)(b_1+b_2)(b_1-b_2)^2 \ when \ N(\ell)=\alpha\\
(\frac{5}{12}\ell^4-\frac{23}{24}\ell^3+\frac{7}{12}\ell^2-\frac{1}{24}\ell)(b_1+b_2)(b_1-b_2)^2 \ when \ N(\ell)=1-\frac{1}{\ell^4}\\
\frac{75}{2}(b_1+b_2)(b_1-b_2)^2 \ when \ \ell=5 \ and \ N(5)=\frac{12}{125}\\
45(b_1+b_2)(b_1-b_2)^2 \ when \ \ell=5 \ and \ N(5)=\frac{96}{625}\\
\frac{155}{2}(b_1+b_2)(b_1-b_2)^2 \ when \ \ell=5 \ and \ N(5)=\frac{252}{625}\\
\end{array} \right.\\
>~&0
\end{align*}
if $b_1 \neq b_2$ and since $2\ell
\binom{\ell}{4}-2\binom{\ell}{3}\ell^2+\binom{\ell}{2}^2\ell=\frac{\ell^2(\ell-1)}{2}>0$
and $\frac{1}{\ell}\geq (b_1+b_2)$.  This completes the proof of Claim
\ref{cl8}.\qed\\
We will apply Claim \ref{cl8} to prove the following claim.

\begin{claim}\label{cl9}
Let $i,\ j,\ 1 \leq i < j \leq q$ be a pair of integers. Let $A(b_1, b_2,\ldots,b_q)$
and $B(b_1, b_2,\ldots,b_q)$ be given as in Claim \ref{cl8}.

Case 1. If $A(b_1, b_2,\ldots,b_q) >0$ then $b_i=b_j$;

Case 2. If $A(b_1, b_2,\ldots,b_q) \leq 0$, then either $b_i=b_j$, or $ min\{b_i,b_j\}=0$.
\end{claim}

The proof of Claim \ref{cl9} (based on Claim \ref{cl8}) can be given by exactly the
same lines as in the proof of Claim 4.5 in \cite{P1} and is omitted here. \qed\\

\emph{Proof of Claim \ref{cl7}.} By Claim \ref{cl9}, either
$b_1=b_2=\cdots=b_q=\frac{1}{\ell q}$ or for some integer $p < q,\
b_{i_1}=b_{i_2}=\cdots=b_{i_p}=\frac{1}{\ell p}$ and other $b_i=0$.

Now we compare $R(\frac{1}{\ell q},\frac{1}{\ell q},\ldots,\frac{1}{\ell
q})=\frac{N(\ell,q)}{120}$ and $R(\frac{1}{\ell p},\frac{1}{\ell p},\ldots,\frac{1}{\ell
p},0,\ldots,0)=\frac{N(\ell,p)}{120}$. It is sufficient to show that $N(\ell,p) \leq
N(\ell,q)$ when $1 \leq p \leq q$. Note that condition (\ref{equa3}) implies that
$N(\ell,1) \leq N(\ell,q)$. Hence it is sufficient to show that $N(\ell,p) \leq
N(\ell,q)$ when $2 \leq p \leq q$ for each of the five choices of $N(\ell)$. In each
case, we view $N(\ell,q)$ as a function with one variable $q$.

{\bf Case a.} $N(\ell)=\alpha$ and $q \geq 2\ell^2 +2\ell$.

In this case, the derivative of $N(\ell,q)$ with respect to $q$ is
\begin{align*}
\frac{d(N(\ell,q))}{dq}&=\frac{10}{\ell q^2}-\frac{70}{\ell^2 q^3}+\frac{150}{\ell^3 q^4}-\frac{16}{\ell^4 q^5}-\frac{40}{\ell q^5}
+\frac{140}{\ell^2 q^5}-\frac{180}{\ell^3 q^5}\\
&=\frac{1}{\ell^4 q^5}(10\ell^3 q^3-70\ell^2 q^2+150\ell q-16-40\ell^3+140\ell^2-180\ell).
\end{align*}
Let $h_1(q)=10\ell^3 q^3-70\ell^2 q^2+150\ell q-16-40\ell^3+140\ell^2-180\ell$, then
$h_1'(q)=30\ell^3 q^2-140\ell^2 q+150\ell,$ $h_1''(q)=60\ell^3 q-140\ell^2.$ Note that
$h_1''(q)>0$ when $q \geq 2,\ \ell \geq 2$, so $h_1'(q)$ increases when $q \geq 2,\ \ell
\geq 2$. By a direct calculation, $h_1'(2)>0$ when $\ell \geq 2$, thus, $h_1(q)$ increases
when $q \geq 2,\ \ell \geq 2$. Since,
$h_1(2)=40\ell^3-140\ell^2+120\ell-16>0$ when $q \geq 2,\ \ell \geq 3$, we know that $N(\ell,q$)
increases when $q \geq 2,\ \ell \geq 3$. When $\ell =2$, by a direct calculation,
$h_1(3)>0$, so $N(2,q$) increases when $q \geq 3$. Also one can easily check that $N(2,2)\leq
N(2,q)$ since $q \geq 2\ell^2 +2\ell$. So $N(\ell,p)\leq N(\ell,q)$ for $2 \leq p \leq
q$.

{\bf Case b.} $N(\ell)=1-\frac{1}{\ell^4}$ and $q \geq 10\ell^3$.

In this case, the derivative of $N(\ell,q)$ with respect to $q$ is
\begin{align*}
\frac{d(N(\ell,q))}{dq}&=\frac{10}{\ell q^2}-\frac{70}{\ell^2 q^3}+\frac{150}{\ell^3 q^4}+\frac{4}{\ell^4 q^5}-\frac{40}{\ell q^5}
+\frac{140}{\ell^2 q^5}-\frac{200}{\ell^3 q^5}\\
&=\frac{1}{\ell^4 q^5}(10\ell^3 q^3-70\ell^2 q^2+150\ell q+4-40\ell^3+140\ell^2-200\ell).
\end{align*}
Let $h_2(q)=10\ell^3 q^3-70\ell^2 q^2+150\ell q+4-40\ell^3+140\ell^2-200\ell$, then
$h_2'(q)=30\ell^3 q^2-140\ell^2 q+150\ell,$ $h_2''(q)=60\ell^3 q-140\ell^2.$ Note that
$h_2''(q)>0$ when $q \geq 2,\ \ell \geq 2$, so $h_2'(q)$ increases when $q \geq 2,\ \ell
\geq 2$. By a direct calculation, $h_2'(2)>0$ when $\ell \geq 2$, thus, $h_2(q)$ increases
when $q \geq 2,\ \ell \geq 2$. Since,
$h_2(2)=40\ell^3-140\ell^2+100\ell+4>0$ when $q \geq 2,\ \ell \geq 3$, we know that $N(\ell,q$)
increases when $q \geq 2,\ \ell \geq 3$. When $\ell =2$, by a direct calculation,
$h_2(3)>0$, so $N(2,q$) increases when $q \geq 3$. Also one can easily check that $N(2,2)\leq N(2,q)$ since $q \geq 10\ell^3$. So $N(\ell,p)\leq N(\ell,q)$ for $2 \leq p \leq q$.

{\bf Case c.} $N(\ell)=\frac{12}{125}$ and $\ell=5$.

In this case, the derivative of $N(5,q)$ with respect to $q$ is
\begin{align*}
\frac{d(N(\ell,q))}{dq}=\frac{2}{q^2}-\frac{14}{5q^3}+\frac{6}{5q^4}-\frac{48}{125q^5}=\frac{1}{125q^5}(250q^3-350q^2+150q-48)\geq
0
\end{align*}
when $q \geq 2$. This proves that $N(5,q)$ increases as $q \geq 2$ increases. So
$N(5,p)\leq N(5,q)$ for $2 \leq p \leq q$.

{\bf Case d.} $N(\ell)=\frac{96}{625}$ and $\ell=5$.

In this case , the derivative of $N(5,q)$ with respect to $q$ is
\begin{align*}
\frac{d(N(\ell,q))}{dq}=\frac{2}{q^2}-\frac{14}{5q^3}+\frac{6}{5q^4}-\frac{384}{625q^5}=\frac{1}{625q^5}(1250q^3-1750q^2+750q-384)\geq
0
\end{align*}
when $q \geq 2$. This proves that $N(5,q)$ increases as $q \geq 2$ increases. So
$N(5,p)\leq N(5,q)$ for $2 \leq p \leq q$.

{\bf Case e.} $N(\ell)=\frac{252}{625}$ and $\ell=5$.

In this case , the derivative of $N(5,q)$ with respect to $q$ is
\begin{align*}
\frac{d(N(\ell,q))}{dq}=\frac{2}{q^2}-\frac{14}{5q^3}+\frac{6}{5q^4}-\frac{1008}{625q^5}=\frac{1}{625q^5}(1250q^3-1750q^2+750q-1008)\geq
0
\end{align*}
when $q \geq 2$. This proves that $N(5,q)$ increases as $q \geq 2$ increases. So
$N(5,p)\leq N(5,q)$ for $2 \leq p \leq q$.

The proof is thus complete.\qed\\

{\bf Acknowledgments.} We also would like to thank the referees for
recommending various improvements of the manuscript. The authors was
supported by NSFC and the ``973" project.

\end{document}